\documentclass[trsc,nonblindrev]{informs3-neutral}

\OneAndAHalfSpacedXI 

\usepackage{booktabs}
\usepackage{longtable}
\usepackage{multirow}
\usepackage{pdflscape}
\usepackage{subfigure}
\usepackage{float}
\usepackage{colortbl}
\usepackage[textsize=scriptsize]{todonotes}
\newcommand{\cO}{{\mathcal{O}}}
\usepackage{nccmath}
\usepackage{enumitem}

\usepackage[boxed,linesnumbered]{algorithm2e}
\DontPrintSemicolon
\SetAlCapSkip{0.5em}
\SetKwComment{tcp}{$\triangleright$\ }{}
\SetKwRepeat{Do}{do}{while}

\SetCommentSty{mycommfont}
\usepackage[normalem]{ulem}

\let\oldnl\nl
\newcommand{\nonl}{\renewcommand{\nl}{\let\nl\oldnl}}

\usepackage{array}
\newcolumntype{H}{>{\setbox0=\hbox\bgroup}c<{\egroup}@{}}
\newcolumntype{C}[1]{>{\centering\arraybackslash}p{#1}}

\usepackage{natbib}
 \bibpunct[, ]{(}{)}{,}{a}{}{,}%
 \def\bibfont{\small}%
 \def\newblock{\ }%

\usepackage{siunitx} 
\TheoremsNumberedThrough     

\EquationsNumberedThrough    

\MANUSCRIPTNO{TS-2021-0280} 
                 
\begin{document}


\RUNAUTHOR{Pacheco et al.}

\RUNTITLE{Exponential-Size Neighborhoods for the PDTSP}

\TITLE{Exponential-Size Neighborhoods for the Pickup-and-Delivery Traveling Salesman Problem}

\ARTICLEAUTHORS{%
\AUTHOR{Toni Pacheco}
\AFF{Departamento de Inform{\'a}tica, Pontif\'{i}cia Universidade Cat\'{o}lica do Rio de Janeiro,\\ \EMAIL{tpacheco@inf.puc-rio.br}}
\AUTHOR{Rafael Martinelli}
\AFF{Departamento de Engenharia Industrial, Pontif\'{i}cia Universidade Cat\'{o}lica do Rio de Janeiro, \\ \EMAIL{martinelli@puc-rio.br}}
\AUTHOR{Anand Subramanian}
\AFF{Departamento de Sistemas de Computa{\c c}{\~a}o, Centro de Inform{\'a}tica, Universidade Federal da Para{\'i}ba, \\ \EMAIL{anand@ci.ufpb.br}}
\AUTHOR{T{\'u}lio A. M. Toffolo}
\AFF{Department of Computing, Federal University of Ouro Preto,\\ \EMAIL{tulio@toffolo.com.br}}
\AUTHOR{Thibaut Vidal}
\AFF{CIRRELT \& SCALE-AI Chair in Data-Driven Supply Chains,\\ D{\'e}partement de Math{\'e}matiques et de G{\'e}nie Industriel, {\'E}cole Polytechnique de Montr{\'e}al,\\
Departamento de Inform{\'a}tica, Pontif\'{i}cia Universidade Cat\'{o}lica do Rio de Janeiro,\\
\EMAIL{thibaut.vidal@cirrelt.ca}}
}

\ABSTRACT{%
Neighborhood search is a cornerstone of state-of-the-art traveling salesman and vehicle routing metaheuristics.
While neighborhood exploration procedures are well developed for problems with individual services, their counterparts for one-to-one pickup-and-delivery problems have been more scarcely studied. A direct extension of classic neighborhoods is often inefficient or complex due to the necessity of jointly considering service pairs. To circumvent these issues, we introduce major improvements to existing neighborhood searches for the pickup-and-delivery traveling salesman problem and new large neighborhoods. We show that the classical \textsc{Relocate-Pair} neighborhood can be fully explored in $\cO(n^2)$ instead of $\cO(n^3)$ time. We adapt the \textsc{4-Opt} and \textsc{Balas-Simonetti} neighborhoods to consider precedence constraints. Moreover, we introduce an exponential-size neighborhood called 2k-\textsc{Opt}, which includes all solutions generated by multiple nested \textsc{2-Opts} and can be searched in $\cO(n^2)$ time using dynamic programming. We conduct extensive computational experiments, highlighting the significant contribution of these new neighborhoods and speed-up strategies within two classical metaheuristics. Notably, our approach permits to repeatedly solve small pickup-and-delivery problem instances to optimality or near-optimality within milliseconds, and therefore it represents a valuable tool for time-critical applications such as meal delivery or mobility on demand.
}

\KEYWORDS{Local search, Large neighborhood search, Dynamic programming, Computational complexity, Pickup-and-delivery problem} 

\maketitle


%
\section{Introduction}
\label{sec:Intro}

Neighborhood search holds a central place in all modern metaheuristics for routing problems \citep{Applegate2006,Vidal2012a,Laporte2014a}, to such an extent that all current state-of-the-art metaheuristics (e.g., \citealt{Nagata2013,Vidal2012b,Subramanian2013,Christiaens2019}) rely on efficient local or large neighborhood search algorithms for achieving high quality solutions. One general thumb rule for success in local-search based metaheuristics is to search \emph{large and fast}. Indeed, larger neighborhoods generally permit further solution refinements, whereas faster exploration techniques permit more restarts (i.e., from initial solutions generated by crossover or perturbation) to diversify the search. Finding a good compromise between neighborhood size and exploration effort is a challenging task, which often involves experimenting with various neighborhoods, pruning rules, and speed-up techniques.

Traveling salesman and vehicle routing problem variants are typically classified in accordance to their service types \citep{Vidal2020}. In one-to-many-to-one \mbox{(1-M-1)} problems (including one-to-many \mbox{1-M} and many-to-one \mbox{M-1} as a special case), each service originates or returns to the depot. Effective local searches in that setting are based on simple exchanges or relocations of customer visits or replacements of arcs. Many neighborhoods for these problems (such as \textsc{2-Opt}, \textsc{Swap} and \textsc{Relocate}) contain $\cO(n^2)$ solutions in their classical form before any other neighborhood restriction. In contrast, in one-to-one pickup-and-delivery problems, every pickup service is paired with its associated delivery, and moves often need joint optimizations of \emph{visit pairs} to remain feasible, leading to neighborhoods which are one order of magnitude more complex (e.g., containing $\cO(n^3)$ solutions) and whose complete exploration can be challenging for large-scale instances.

In this article, we introduce new neighborhood searches for one-to-one pickup-and-delivery problems. More precisely, we focus on the pickup-and-delivery traveling salesman problem (PDTSP), which is the canonical problem in this family, and make the following contributions:
\begin{itemize}
\item We design an efficient evaluation strategy which permits to explore the \textsc{Relocate Pair} neighborhood in $\cO(n^2)$ instead of $\cO(n^3)$ time.
\item We introduce a neighborhood of exponential size called 2k-\textsc{Opt}. This neighborhood contains all solutions formed with \emph{nested} \textsc{2-Opts} and can be explored in $\cO(n^2)$ time. This neighborhood overcomes the natural tendency of \textsc{2-Opt} to create infeasible solutions due to route reversals and precedence constraints in pickup-and-delivery problems.
\item We adapt the restricted \textsc{4-Opt} neighborhood of \cite{Glover1996a} and the neighborhood of \cite{Balas2001} to the PDTSP.
\item We integrate all these neighborhoods within a sophisticated local search and evolutionary metaheuristic: the hybrid genetic search (HGS) of \cite{Vidal2012,Vidal2012b}. As visible in our computational experiments on the PDTSP, the proposed solution approach performs remarkably well on the 163 classical instances from \cite{Renaud2000a} and \cite{Dumitrescu2010}, identifying better or equal solutions than all previous algorithms in all cases. We also conduct sensitivity analyses on the individual contributions of each neighborhood. Based on these experiments, the efficient \textsc{Relocate Pair}, the \textsc{Or-Opt} and the adapted \textsc{4-Opt} are highlighted as the most important search operators for this problem.
\end{itemize}

\section{Related Literature}
\label{literature}

\textbf{Neighborhood search.}
Let $X$ be the set of all feasible solutions of a given optimization problem. Formally, a neighborhood is a mapping $N: X \rightarrow 2^X$ associating to each solution $x$ a set of neighbors $N(x) \subseteq X$. Most vehicle routing neighborhoods are built on the definition of a \emph{move} operator (e.g., a \textsc{Relocate}, which consists in changing the position of a service in the routing plan) and contain all solutions reachable from $X$ by a single move. Based on this definition, a local-search procedure iteratively explores the neighborhood of an incumbent solution to find an improved one, replaces the incumbent solution with this solution, and repeats this process until no further improvement exists. The final solution is called a local optimum.

Local searches for 1-M-1 problems typically rely on \textsc{$k$-opt}, \textsc{Relocate} and \textsc{Swap} moves and some of their extensions. A \textsc{$k$-opt} move consists of replacing up to $k$ arcs in the solution. \textsc{Relocate} consists of moving a customer visit to another position, whereas \textsc{Swap} exchanges the positions of two customer visits. There are $\cO(k!n^k)$ possible \textsc{$k$-opt} and $\cO(n^2)$ possible \textsc{Relocate}~and~\textsc{Swap}~moves. In contrast, 1-1 pickup-and-delivery problems involve paired services. Neighborhood searches are more complex in this setting since the locations of the pickups and of their corresponding deliveries need to be jointly optimized. Single-visit \textsc{Relocate} moves, for example, often generate infeasible solutions violating precedence or capacity constraints. As a consequence, it is common to explore the \textsc{Relocate-Pair} neighborhood, which consists of jointly relocating a pickup and its associated delivery. However, an enumerative exploration of this neighborhood takes $\cO(n^3)$ time.

Various speed-up techniques can be used to improve neighborhood search. The well-known \cite{Lin1973} algorithm and its implementation by \cite{Helsgaun2000, Helsgaun2009} combines filtering mechanisms and bounds to perform an effective search of \textsc{$k$-opt} of order up to $k=5$ for the traveling salesman problem. Other simple neighborhood restriction strategies consist in restricting the moves to ``close'' locations \citep{Johnson1997,Toth2003}. In a different fashion, \cite{Taillard2019} firstly generate a pool of initial solutions and restrict the subsequent search to the edges belonging at least to one of these initial solutions. Efficient search strategies and spatial data structures may also be exploited to achieve significant speedups \citep{Bentley1992,Glover1996a,DeBerg2020,Lancia2019}. 

Some families of neighborhoods qualified as \emph{very large} contain an exponential number of solutions but can be efficiently searched \citep{Ahuja2002}. This is often done by exploiting combinatorial subproblems, e.g., shortest paths, assignment or network flows \citep{Franceschi2006,Hintsch2018,Vidal2017b,Capua2018} and dynamic programming techniques \citep{Deineko2000,Balas2001,Congram2002}. \cite{Balas2001} exploit the availability of a linear-time dynamic programming algorithm for a restricted variant of the traveling salesman problem (TSP) to design a very large neighborhood. Given an incumbent solution represented as a permutation $\sigma$ and a fixed value $k$, this neighborhood contains all permutations $\pi$ of $\sigma$ such that $\pi$ fulfills $\pi(1) = 1$ and $\pi(i) \leq \pi(j)$ for all $i,j \in \{1,\dots,n\}$ such that $\sigma(i)+k \leq \sigma(j)$. In other words, if $j$ is located at least $k$ positions after $i$ in the original permutation $\sigma$, then $j$ must appear after $i$ in permutation $\pi$.
This neighborhood was later exploited for several variants of vehicle routing and arc routing problems \citep{Irnich2008,Vidal2017b,Gschwind2019,Toffolo2019}.
Finally, ejection chain neighborhoods consist of exploring chained customer relocations. Some exponential-size neighborhoods of this class can be reduced to the search of a negative-cost cycle or path in an auxiliary graph \citep{Glover2006,Vidal2017b}. It is also noteworthy that efficient search procedures for exponential-size neighborhoods exist for the TSP, but that they are theoretically impossible for  other problems such as the quadratic assignment problem unless P=NP (see \citealt{Deineko2000}).\\

\noindent
\textbf{The pickup-and-delivery TSP.}
The PDTSP is the canonical form of a single-vehicle one-to-one pickup-and-delivery problem. This problem is highly relevant in robotics, production and cutting as well as transportation logistics. It has become the focus of renewed research due to its direct applications to mobility-on-demand and meal-delivery services \citep[see, e.g.,][]{ONeil2018a,Yildiz2019}. Formally, the PDTSP can be defined over a graph $G=(V,E)$ in which the vertices $V = 0 \cup P \cup D$ represent service locations and the edges in $E$ indicate possible trips between locations. Let $P=\{1,\ldots,n\}$ and $D=\{n+1,\ldots,2n\}$ represent pickup and delivery vertices, respectively. There are $n$ pickup-and-delivery requests $(i,n+i)$ for $i \in \{1,\ldots,n\}$, each one associated to a pickup location $i \in P$ and a delivery location $n+i\in D$. Let $N = 2n + 1$ be the number of visits in the problem. Vertex $0$ represents an origin location, from which each vehicle should start and return. Finally, each edge $(i,j) \in E$ is associated to a travel cost $c_{ij}$. The goal of the PDTSP is to find a tour of minimum total distance, starting and ending at $0$, in such a way that each pickup and delivery location is visited once, and that each pickup precedes its corresponding delivery in the tour.

The PDTSP is formulated as a binary integer program in Equations \eqref{eq:form} to \eqref{eq-subj6}. This model and notations originates from \cite{Ruland1994} and \cite{Dumitrescu2010}. The binary decision variable $x_{ij}$ takes value $1$ if and only if the vehicle travels on edge $(i,j)\in E$. For each subset of vertices $S\subseteq V$, we define the cut set $\delta(S)=\{(i,j)\in E : i \in S, j \notin S ~\text{or}~ i \notin S, j \in S\}$. Moreover, $\delta(i)$ is used instead of $\delta(\{i\})$ for unit sets and $\displaystyle x(E') = \sum_{(i,j)\in E'} x_{ij}$ for any subset $E' \subseteq E$.
\begin{align}
\text{minimize} \hspace{0.5cm}  & \sum_{(i,j)\in E} c_{ij}x_{ij} \label{eq:form} \\ 
\text{subject to} \hspace{0.5cm} & x_{0,2n+1} = 1 \label{eq-subj2}\\
    & x(\delta(i)) = 2     & \forall\, i \in V \label{eq-subj3}\\
    & x(\delta(S)) \geq 2  & \forall\, S \subseteq V, 3 \leq |S| \leq |V|/2 \label{eq-subj4}\\
    & x(\delta(S)) \geq 4  & \forall\, S \in \mathcal{U} \label{eq-subj5}\\
    & x_{ij} \in \{0,1\}   & \forall\, (i,j) \in E. \label{eq-subj6}
\end{align}

In this formulation, $\mathcal{U}$ is the collection of all subsets $S \subset V$ such that $3 \leq |S| \leq |V|-2$, $0 \in S$, $2n+1 \notin S$, and there exists $i \in P$ such that $i \notin S$ and $n+i \in S$. Objective~\eqref{eq:form} minimizes the total travel distance. Constraint~\eqref{eq-subj2}  permits to model the solution as a flow circulation. Constraints~\eqref{eq-subj3} and \eqref{eq-subj4} are the standard degree and subtour elimination constraints, respectively. Finally, Constraints~\eqref{eq-subj5} enforce every pickup location $i\in P$ to be visited before its corresponding delivery location $n+i$.

Despite the availability of sophisticated integer programmings solvers, an exact solution of the PDTSP remains impractical beyond a few dozens of pickup-and-delivery requests \citep{Dumitrescu2010}. This is in stark contrast with the classical TSP, for which instances with thousands of nodes can be routinely solved to proven optimally \citep{Applegate2009}. As a consequence, metaheuristics remain the most suitable alternative to solve larger instances. \cite{Psaraftis1983,Savelsbergh1990} and \cite{Healy1995} were among the first to study construction and local improvement heuristics for this problem. \cite{Renaud2000a} introduced a construction and local search algorithm, which was later refined by adding seven different perturbation strategies to escape local optima in \cite{Renaud2002a}. Later on, \cite{Dumitrescu2010} used a large neighborhood search (LNS) to generate initial upper bounds for their branch-and-cut algorithm. \cite{Veenstra2017} designed another LNS with five different destruction operators to address the PDTSP and its generalization with handling costs. To our knowledge, this method is regarded as the current state-of-the-art for the PDTSP.

Finally, several PDTSP variants with side constraints have been introduced to accommodate diverse requirements of practical applications, e.g., time windows, capacity constraints, handling costs, and loading constraints. These variants are thoroughly surveyed in \cite{Parragh2008a,Parragh2008}. LIFO and FIFO loading constraints, in particular, have been the focus of several studies \citep{Carrabs2007,Erdogan2009,Cordeau2012a,Pollaris2015}. From a methodological standpoint, LIFO or FIFO can restrict the number of feasible solutions, but also enhance neighborhood searches as they uniquely determine a delivery position within the tour from its pickup position. Overall, despite decades of research, significant progress remains needed regarding efficient neighborhood exploration procedures for the unconstrained PDTSP.

\section{Neighborhood Exploration Methodology} 
\label{sec:neighborhoods}

We propose to rely on six neighborhoods: \textsc{Relocate Pair}, \textsc{2-Opt}, \textsc{Or-Opt}, \textsc{2$k$-Opt}, \textsc{4-Opt}, and \textsc{Balas-Simonetti}. We divide these neighborhoods into two groups according to the way they operate. The first group contains the first three neighborhoods. A common characteristic of these neighborhoods is that it is possible to decompose their exploration, e.g., per pickup-and-delivery pair. Therefore, improving moves can be applied as soon as they are detected. These three neighborhoods are well known, yet we propose a new exploration strategy for \textsc{Relocate Pair} which takes $\mathcal{O}(n^2)$ instead of $\mathcal{O}(n^3)$ time. The second group contains the last three neighborhoods, with exploration procedures based on dynamic programming (DP). A complete DP is typically needed to identify the best move. These neighborhoods are either new (\textsc{2$k$-Opt}), or adaptations of TSP neighborhoods (\textsc{4-Opt} and \textsc{Balas-Simonetti}). The \textsc{2$k$-Opt} neighborhood consists of combining nested \textsc{2-Opt} nested moves in $\mathcal{O}(n^2)$ operations, while the adaptations of \textsc{4-Opt} and \textsc{Balas-Simonetti} neighborhoods allow to evaluate a large set of feasible PDTSP moves in $\mathcal{O}(n^2)$ and $\mathcal{O}(k^2 2^{k-2} n)$ time, respectively. 

All these neighborhoods are restricted to feasible solutions satisfying the precedence constraints between pickups and deliveries. We describe each neighborhood search algorithm in a dedicated subsection. Each search procedure is applied on a feasible incumbent solution $\sigma$, represented as a visit sequence $\sigma$ such that $\sigma(0)=\sigma(2n+1)=0$ and $\sigma(i) \in \{1,\dots,2n\}$ for all $i \in \{1,\dots,2n\}$. In this text, we will also refer to $\sigma_{[i,j]}$ as the subsequence of consecutive visits $(\sigma(i), \ldots, \sigma(j))$.

\subsection{\textsc{Relocate Pair} Neighborhood} 
\label{relocate-pair}

\begin{figure}[htbp]
\centering
\boxed{
\begin{minipage}{0.47\linewidth}
\centering
\small
\vspace*{0.38cm}
\includegraphics[width=\linewidth]{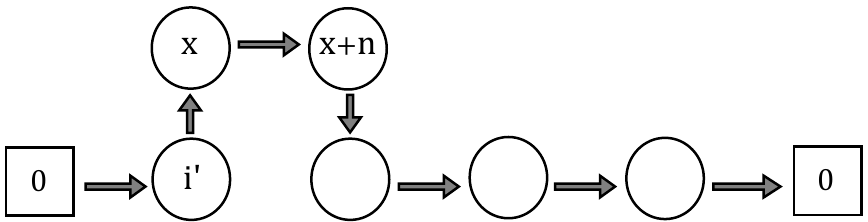}
A) Insertion in Consecutive Positions
\end{minipage}
}
\hspace*{0.2cm}
\boxed{
\begin{minipage}{0.47\linewidth}
\centering
\small
\vspace*{-0.1cm}
\includegraphics[width=\linewidth]{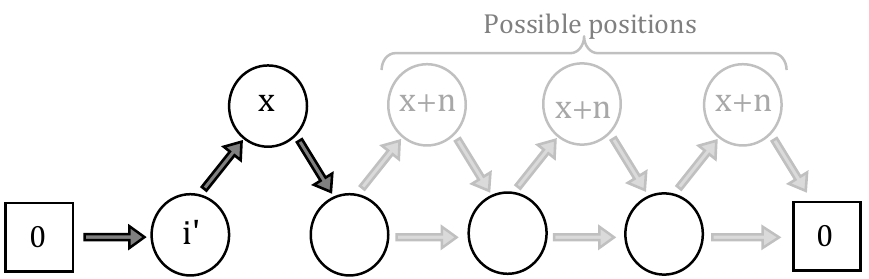}
B) Insertion in Non-Consecutive Positions
\end{minipage}
}
\vspace*{0.15cm}
	\caption{Illustration of the \textsc{Relocate Pair} neighborhood}
	\label{fig:relocate-pair-search}
\end{figure}

The \textsc{Relocate Pair} neighborhood is defined by a \emph{move} operation which consists in relocating a pickup and delivery pair $(x,x+n)$ from its current positions $i$ and $j$ to new positions $i'$ and $j'$ such that $i' < j'$. Our search procedure iterates over the $n$ pickup and delivery pairs $(x,n+x)$ in random order.
For each such pair, it firstly calculates in $\cO(1)$ operations the removal cost:
\begin{equation}
\Delta_\textsc{rem} = 
\begin{cases}
c_{\sigma(i-1),\sigma(j+1)} - c_{\sigma(i-1),\sigma(i)} - c_{\sigma(i),\sigma(j)} - c_{\sigma(j),\sigma(j+1)} & \text{if}~j=i+1\\
c_{\sigma(i-1),\sigma(i+1)} - c_{\sigma(i-1),\sigma(i)} - c_{\sigma(i),\sigma(i+1)} + c_{\sigma(j-1),\sigma(j+1)} - c_{\sigma(j-1),\sigma(j)} - c_{\sigma(j),\sigma(j+1)}& \text{otherwise.}
\end{cases}
\end{equation}
Then, it searches for the best insertion positions $i'$ and $j'$ within the route $\pi$ after removal (such that $\pi(0) = \pi(2n-1) = 0$, and $\pi(i)$ for $i \in \{1,\dots,2n-2\}$ represent customer visits). As illustrated in Figure~\ref{fig:relocate-pair-search}, two possibilities need to be considered.

\noindent
\textbf{Case A)}. The best insertion cost in \textbf{consecutive} positions $\Delta^\textsc{c}_\textsc{add}$ can be found in $\cO(n)$ time by simple inspection, using the following equation:
\begin{equation}
\Delta^\textsc{c}_\textsc{add} = \min_{i' \in \{0,\dots,2n-2\}} \{ c_{\pi(i'),x} + c_{x,n+x} + c_{n+x,\pi(i'+1)} - c_{\pi(i'),\pi(i'+1)} \}.\label{eq:rp:1}
\end{equation}

\noindent
\textbf{Case B)}. The best insertion cost in \textbf{non-consecutive} positions $\Delta^\textsc{nc}_\textsc{add}$ is calculated as:
\begin{equation}
\Delta^\textsc{nc}_\textsc{add} = \hspace*{-0.1cm}  \min_{i' \in \{0,\dots,2n-3\}} \left\{ c_{\pi(i'),x} + c_{x,\pi(i'+1)} - c_{\pi(i'),\pi(i'+1)} + \hspace*{-0.1cm}  \min_{j' \in \{i'+1,\dots,2n-2\}} \left\{
c_{\pi(j'),n+x} + c_{n+x,\pi(j'+1)} - c_{\pi(j'),\pi(j'+1)} \right\} \hspace*{-0.1cm} \right\} \hspace*{-0.05cm}. \vspace*{0.2cm}
\end{equation}
A direct application of this equation takes $\cO(n^2)$ time \emph{for each pair} $(x,x+n)$, leading to a complete neighborhood evaluation in $\cO(n^3)$ time. To avoid this high computational effort, we observe that: \vspace*{0.2cm}
\begin{equation}
\Delta^\textsc{nc}_\textsc{add} = \min_{i'} \left\{ \Delta(i') + \Phi(i') \right\},
\text{ with }
\begin{cases}
\Delta(i') = c_{\pi(i'),x} + c_{x,\pi(i'+1)} - c_{\pi(i'),\pi(i'+1)} \\
\Phi(i') = \smash{\min\limits_{j' \in \{i'+1,\dots,2n-2\}}} \left\{ c_{\pi(j'),n+x} + c_{n+x,\pi(j'+1)} - c_{\pi(j'),\pi(j'+1)} \right\}, \vspace*{0.2cm}
\end{cases}
\end{equation}
and that all the values $\Phi(i')$ can be calculated by backward recursion, in $\cO(n)$ time overall: \vspace*{0.2cm}
\begin{equation}
\Phi(i') =  
\begin{cases}
c_{\pi(2n-2),n+x} + c_{n+x,0} - c_{\pi(2n-2),0} & \text{if } i' = 2n-3\\
\min \{ \Phi(i'+1), \smash{c_{\pi(i'+1),n+x} + c_{n+x,\pi(i'+2)} - c_{\pi(i'+1),\pi(i'+2)}} \} & \text{otherwise.} \vspace*{0.2cm}
\end{cases} \label{eq:rp:2}
\end{equation}

As a consequence $\Delta^\textsc{c}_\textsc{add}$ and $\Delta^\textsc{nc}_\textsc{add}$ are calculated in $\cO(n)$ time. The best \textsc{Relocate Pair} move for the pickup and delivery pair $(x,x+n)$ has a cost of $\Delta = \Delta_\textsc{rem} + \min \{\Delta^\textsc{c}_\textsc{add}, \Delta^\textsc{nc}_\textsc{add} \}$ and is applied if $\Delta < 0$. Overall, running this search procedure for all pickup and delivery pairs requires $\cO(n \times n) = \cO(n^2)$ time and $\cO(n)$ space.

\subsection{\textsc{2-Opt} Neighborhood} 
\label{sec:2opt}

The \textsc{2-Opt} neighborhood is defined by a move which replaces two edges with two new ones.
Given the current solution $\sigma$ and two positions $i$ and $j$ such that $i + 2 < j$, a \textsc{2-Opt} move on these positions reverses the visit sequence $\sigma_{[i+1,j-1]}$ and replaces arcs \mbox{$(\sigma(i)$, $\sigma(i+1))$} and \mbox{$(\sigma(j-1)$, $\sigma(j))$} by arcs \mbox{$(\sigma(i),\sigma(j-1))$} and \mbox{$(\sigma(i+1), \sigma(j))$}. This leads to a solution improvement when the following cost difference is negative:
\begin{equation}
\Delta_\textsc{2opt}(i,j) = c_{\sigma(i),\sigma(j-1)} + c_{\sigma(i+1),\sigma(j)} - c_{\sigma(i),\sigma(i+1)} - c_{\sigma(j-1),\sigma(j)}. \label{eq:cost2opt}
\end{equation}

Evaluating all \textsc{2-Opts} takes $\mathcal{O}(n^2)$ time.
However, since any \textsc{2-Opt} for positions $i$ and $j$ reverses the visit sequence $\sigma_{[i+1,j-1]}$, it leads to an infeasible solution whenever there exists a pickup and its associated delivery within the subsequence. These infeasible cases can be efficiently filtered out during the search when iterating over $i$ and $j$, by stopping the inner loop whenever position $j$ corresponds to a delivery and the associated pickup belongs to $\sigma_{[i+1,j-1]}$.

\subsection{\textsc{Or-Opt} Neighborhood}
\label{sec:oropt}

Each \textsc{Or-Opt} move \citep{Or1976} consists in relocating a sequence of visits $\sigma_{[i,j]}$ at a different position~$k$. There are $\mathcal{O}(n^3)$ such moves. Since each \textsc{Or-Opt} move breaks three edges at most, \textsc{Or-Opt} is effectively a subset of \textsc{3-Opt}.

To cut down the computational complexity of neighborhood exploration, we opt to restrict this neighborhood to sequences $\sigma_{[i,j]}$ no longer than a constant $k_\textsc{or}$, and we allow joint sequence relocations and reversals. This restriction is meaningful since relocations of long sequences tend to be frequently infeasible. Moreover, with this limitation, the evaluation of all moves can be done by simple enumeration~in~$\mathcal{O}(k_\textsc{or} n^2)$ time.

\subsection{2k-\textsc{Opt} Neighborhood} 
\label{sec:2k-opt}

\begin{figure}[htbp]
\centering
\boxed{
\includegraphics[width=0.62\linewidth]{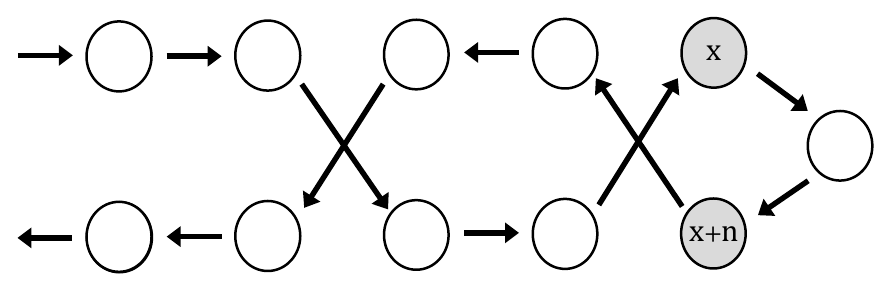}
}
\caption{Example of a setting in which a 2k-\textsc{Opt} move is essential to improve distance.}
\label{fig:k2opt-extramoves}
\end{figure}

The \textsc{2-Opt} neighborhood is widely used for the classical TSP, but its efficiency for the PDTSP is drastically limited since it reverses a visit sequence and therefore often violates precedence constraints. For example, Figure~\ref{fig:k2opt-extramoves} represents a fragment of solution (route). In this fragment, all pickups (respectively, deliveries) are paired with deliveries (respectively, pickups) located outside of the sequence, except for one pickup-and-delivery pair denoted as $x$ and $x+n$. Assuming that distances are Euclidean, the two crossings induced by the current solution visibly lead to extra distance. Moreover, there exists two \textsc{2-Opt} moves capable of uncrossing the solution, but each of these moves, in isolation, would reverse the order of $x$ and $x+n$ and render the solution infeasible due to precedence constraints. Obviously, the solution, in such a case, is to perform both \textsc{2-Opt} moves simultaneously. As seen in the following, this is a special case of the \textsc{2k-Opt} neighborhood, and we propose dynamic programming algorithms to search this neighborhood efficiently.

We will demonstrate that the best 2k-\textsc{Opt} move can be found by dynamic programming in $\mathcal{O}(n^2)$ operations.
For any pair of indices $(i,j)$ such that $1 \leq i \leq j \leq 2n$, let $F(i,j)$ be the cost of the best combination of \textsc{2-Opt} in the visit sequence $\sigma_{[i,j]}$. In a similar fashion, let $R(i,j)$ be the cost of the best combination of \textsc{2-Opt} within the reversed sequence~$\pi_{[i,j]} = (\sigma(j)\dots,\sigma(i))$.
Evaluating $F(i,j)$ can be done by induction, selecting the best out of three cases:
\begin{enumerate}
    \item[1)] Searching for the best \textsc{2-Opt} combination in $\sigma_{[i+1,j]}$; 
    \item[2)] Searching for the best \textsc{2-Opt} combination in $\sigma_{[i,j-1]}$;
    \item[3)] Applying \textsc{2-Opt}$(i,j)$ and searching for additional moves in the reversed sequence $\pi_{[i+1,j-1]}$.
\end{enumerate}
This policy can be represented by the following equations, in which $\Delta_\textsc{2opt}(i,j)$ represents the cost of a \textsc{2-Opt} between $i$ and $j$ as defined in Equation~(\ref{eq:cost2opt}).
\begin{fleqn}
\def\arraystretch{1.7}
    \begin{align}
        ~&F(i,j) = \left\{ \begin{array}{ll} 
            0, & \text{if}~j\leq i+1\\
            \min \left\{ \begin{array}{l}
                F(i+1,j),\\~F(i,j-1),\\~~\Delta_\textsc{2opt}(i,j) + R(i+1,j-1)
            \end{array} \right\}, &\text{otherwise}
        \end{array} \right.
        \label{eq:recurence_f}
    \end{align}
\end{fleqn}
\vspace*{0.3cm}

\noindent
In a similar fashion, the best policy for $R(i,j)$ can be calculated as:
\begin{fleqn}
\def\arraystretch{1.7}
    \begin{align}
        ~&
        R(i,j) = \left\{ \begin{array}{ll} 
             \infty, & \begin{array}{l} \text{if}~\sigma(i)+n \in \sigma_{[i+1,j]}~\text{or}~\sigma(j)-n \in \sigma_{[i,j-1]}\end{array}\\
            0, & \begin{array}{l} \text{if}~j = i \end{array}\\ 
            0, & \begin{array}{l} \text{if}~j = i+1~\text{and}~\sigma(j)\neq n+\sigma(i) \end{array}\\ 
            \infty, & \begin{array}{l} \text{if}~j = i+1~\text{and}~\sigma(j)=n+\sigma(i) \end{array}\\ 
            \min \left\{ \begin{array}{l}
                R(i+1,j),\\~R(i,j-1),\\~~\Delta_\textsc{2opt}(i,j)+F(i+1,j-1)
            \end{array} \right\}, & \text{otherwise.}
        \end{array} \right. 
        \label{eq:recurence_r}
    \end{align}
\end{fleqn}
In this recursion, the first statement eliminates cases in which the node $\sigma(i)$ is paired with a delivery in $\sigma_{[i+1,j]}$, or in which $\sigma(j)$ is paired with a pickup in $\sigma_{[i,j-1]}$. Indeed, regardless of the internal \textsc{2-Opt} choices, $\sigma(j)$ and $\sigma(i)$ are guaranteed to be respectively the first and the last visited. The next three statements cover the cases of sequences containing one or two visits, and the final statement performs the recursive calls in a similar fashion as Equation~(\ref{eq:recurence_f}).

The number of possible states of $F(i,j)$ and $R(i,j)$ is in $\cO(n^2)$ and the evaluation of the best cost for each state is done in constant time through Equations (\ref{eq:recurence_f}--\ref{eq:recurence_r}). Overall, a direct application of these equations gives an algorithm to find the best 2k-\textsc{Opt} move in $\cO(n^2)$ time and space.

\subsection{\textsc{4-Opt} Neighborhood}

Moves of a higher order such as \textsc{4-Opt} have the potential to lead to structurally different solutions that are not attainable with simpler moves. However, a complete exploration of this neighborhood is prohibitively long. Some previous works have opted to explore only a small fraction of this neighborhood. \cite{Renaud2000a,Renaud2002a}, for example, explored a subset of $\cO(n^2)$ moves.

\begin{figure}[htbp]	
	\centering
	\subfigure[Type 1: Two intertwined \textsc{2-Opt-d}]
    {\includegraphics[width=0.4\linewidth]{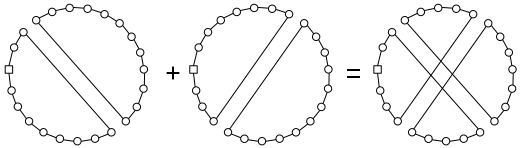}
	\label{fig:4opt-type1}}\hspace{10mm}
    \subfigure[Type 2: Intertwined \textsc{2-Opt-c} and \textsc{2-Opt-d}]
    {\includegraphics[width=0.4\linewidth]{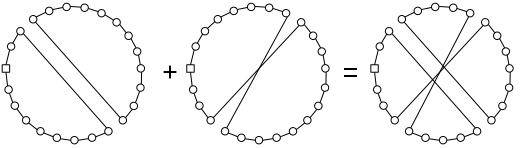}
	\label{fig:4opt-type2}}
    \subfigure[Type 3: Two intertwined \textsc{2-Opt-c}]
    {\includegraphics[width=0.4\linewidth]{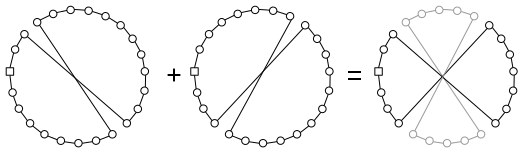}
	\label{fig:4opt-type3}}\hspace{15mm} 
    \subfigure[Type 4: Two disjoint \textsc{2-Opt-d}]
    {\includegraphics[width=0.4\linewidth]{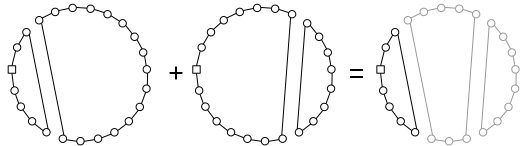}
	\label{fig:4opt-type4}} 
    \subfigure[Type 5: Disjoint \textsc{2-Opt-c} and \textsc{2-Opt-d}]{\includegraphics[width=0.4\linewidth]
    {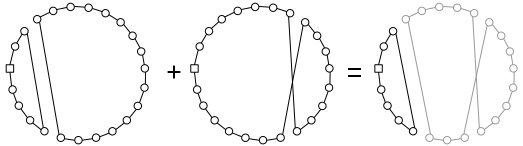}
	\label{fig:4opt-type5}}\hspace{15mm} 
	\subfigure[Type 6: Two disjoint \textsc{2-Opt-c}]
    {\includegraphics[width=0.4\linewidth]{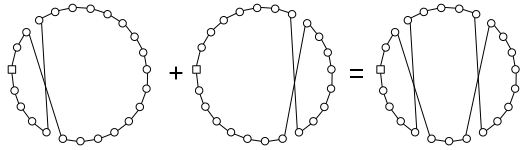}
	\label{fig:4opt-type6}}
    \vspace{2.5mm} 
	\caption{4-Opt moves obtained as a combination of \textsc{2-Opt-c} and \textsc{2-Opt-d}.\label{fig:4opt-types}} 
\end{figure}

We take a fairly different approach to explore this neighborhood, and instead rely on an extension of the search strategy of \cite{Glover1996a} originally designed for the TSP. As discussed in this section, this allows exploring a subset of the \textsc{4-Opt} neighborhood with up to $\Theta(n^4)$ moves in $\cO(n^2)$ time. This search strategy can be described as a search for an improving alternating cycle in an auxiliary graph. Here we will adopt a direct presentation in a recursive form. We focus on the \textsc{4-Opt} moves that are obtained as a combination of a pair of \emph{connecting} or \emph{disconnecting} \textsc{2-Opt} moves, referred as \textsc{2-Opt-c} and \textsc{2-Opt-d} throughout this section. As listed in Figure~\ref{fig:4opt-types}, six possible configurations exist.

Only three of these configurations (Types 1, 2, and 6) produce connected solutions:
\begin{itemize}[nosep]
\item Type-1 moves are obtained by intertwining two \textsc{2-Opt-d} moves. Such moves are commonly used as a perturbation operator and called ``double bridge'' in the TSP and vehicle routing literature \citep{Subramanian2013}.
\item Type-2 moves are obtained by intertwining a \textsc{2-Opt-c} move with a \textsc{2-Opt-d} move. This type will be divided into two cases (Types 2A and 2B) depending on whether \textsc{2-Opt-d} is applied first or second. Some visit sequences are reversed when applying this move.
\item Type-6 moves consist of two disjoint \textsc{2-Opt-c} moves, and represent a special case of the \textsc{2k-Opt} neighborhood described in Section~\ref{sec:2k-opt}. For this reason, we focus exclusively on Type-1 and Type-2 moves. 
\end{itemize}

\begin{figure}[htbp]	
    \centering
    \begin{minipage}{0.3\textwidth}
    \centering
    \includegraphics[width=\textwidth]{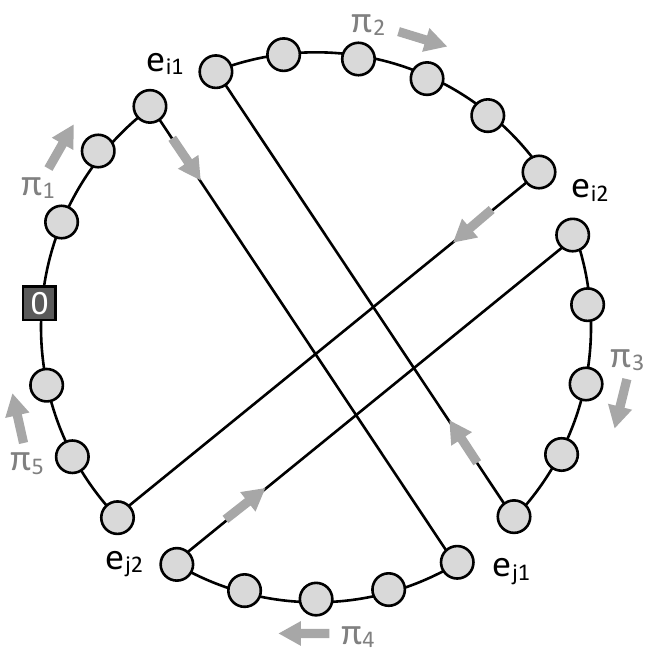}

    \textbf{Type 1}
    \end{minipage}
    \hfill
    \begin{minipage}{0.3\textwidth}
    \centering
    \includegraphics[width=\textwidth]{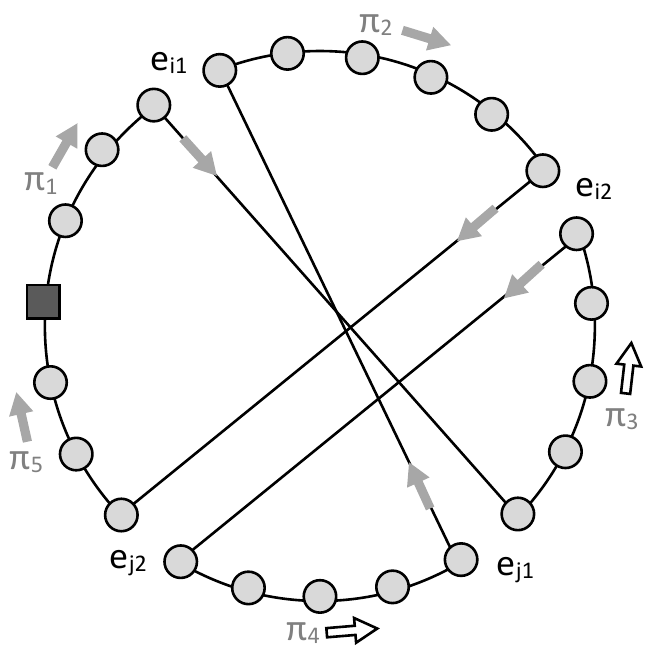}

    \textbf{Type 2A}
    \end{minipage}
    \hfill
    \begin{minipage}{0.3\textwidth}
    \centering
    \includegraphics[width=\textwidth]{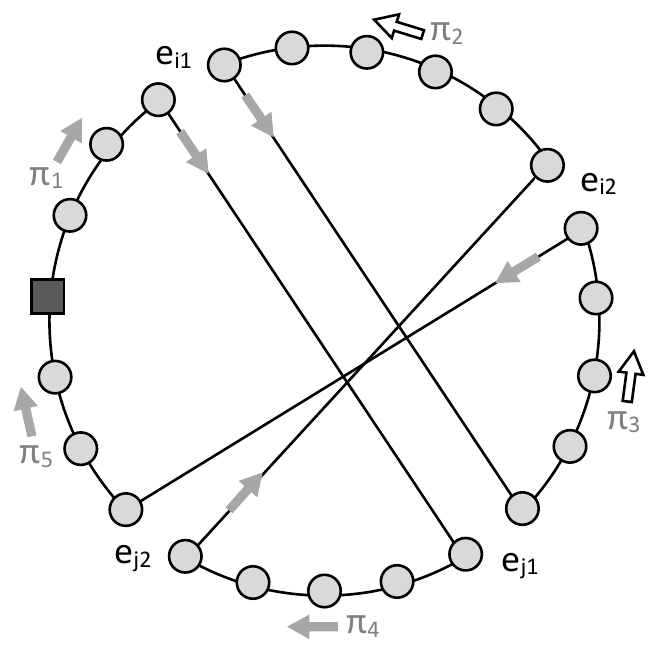}

    \textbf{Type 2B}
    \end{minipage}
    \vspace*{0.5cm}
    \caption{Different reconnections and sequence orders corresponding to different types of \textsc{4-opt} moves\label{configs-reconnection}}
\end{figure}

Each \textsc{4-Opt} of a given type is characterized by four indices $0 \leq {i_1}<{i_2}<{j_1}<{j_2} \leq 2n$ representing the removal of edges \mbox{$e_{i_1} = (\sigma(i_1), \sigma(i_1+1))$}, \mbox{$e_{i_2} = (\sigma(i_2), \sigma(i_2+1))$}, \mbox{$e_{j_1}= (\sigma(j_1), \sigma(j_1+1))$} and \mbox{$e_{j_2} = (\sigma(j_2), \sigma(j_2+1))$}.
After removal, the incumbent solution $\sigma$ is split into five visit sequences: $\pi_1 = \sigma_{[0,i_1]}$, $\pi_2 = \sigma_{[i_1+1,i_2]}$, $\pi_3 = \sigma_{[i_2+1,j_1]}$, $\pi_4 = \sigma_{[j_1+1,j_2]}$ and $\pi_5 = \sigma_{[j_2+1,2n]}$. During reconstruction, the extremities of edge $e_{i_1}$ are reconnected with those of edge $e_{j_1}$. There are two possible ways to proceed: $\textsc{2-Opt-c}(i_1,j_1)$ creates edges \mbox{$(\sigma(i_1), \sigma(j_1))$} and \mbox{$(\sigma(i_1+1), \sigma(j_1+1))$}, whereas $\textsc{2-Opt-d}(i_1,j_1)$ creates edges \mbox{$(\sigma(i_1), \sigma(j_1+1))$} and \mbox{$(\sigma(i_1+1), \sigma(j_1))$}. Likewise, the extremities of edges $e_{i_2}$ and $e_{j_2}$ can be reconnected in two possible ways. This leads to the following sequence orientations and routes for each move type (Figure~\ref{configs-reconnection}):
\begin{itemize}
    \item Type 1: $\textsc{2-Opt-d}(i_1,j_1)$ and $\textsc{2-Opt-d}(i_2,j_2)$, producing route $\pi_1 \oplus \pi_4 \oplus \pi_3 \oplus \pi_2 \oplus \pi_5$;
	\item Type 2A:  $\textsc{2-Opt-c}(i_1,j_1)$ and $\textsc{2-Opt-d}(i_2,j_2)$, producing route $\pi_1 \oplus \textrm{r}(\pi_{3}) \oplus \textrm{r}(\pi_{4}) \oplus \pi_2 \oplus \pi_5$;
    \item Type 2B:  $\textsc{2-Opt-d}(i_1,j_1)$ and $\textsc{2-Opt-c}(i_2,j_2)$, producing route $\pi_1 \oplus \pi_4 \oplus \textrm{r}(\pi_{2}) \oplus \textrm{r}(\pi_{3}) \oplus \pi_5$;
\end{itemize}
where \textrm{r}() represents the reversal of a sequence and $\oplus$ is the operation of concatenation of two sequences.

\noindent
The cost of each $\textsc{2-Opt-d}$ and $\textsc{2-Opt-c}$ move can be calculated as:
\begin{align}
\Delta^\textsc{d}_\textsc{2opt}(i,j) &= c_{\sigma(i),\sigma(j+1)} + c_{\sigma(i+1),\sigma(j)} - c_{\sigma(i),\sigma(i+1)} - c_{\sigma(j),\sigma(j+1)}.  \\
\Delta^\textsc{c}_\textsc{2opt}(i,j) &= c_{\sigma(i),\sigma(j)} + c_{\sigma(i+1),\sigma(j+1)} - c_{\sigma(i),\sigma(i+1)} - c_{\sigma(j),\sigma(j+1)},
\end{align}
and the cost of the best \textsc{4-opt} move for each fixed edge pair $i_2$ and $j_2$ can be calculated as:
\begin{align}
\Delta^\textsc{Type-1}_\textsc{4opt}(i_2,j_2) &= \Delta^\textsc{D}_\textsc{2opt}(i_2,j_2) + \min_{i_1 \in \{0,\dots,i_2-1\}} \left\{ \min_{j_1 \in \{i_2+1,\dots,j_2-1\}} \left\{ \Delta^\textsc{D}_\textsc{2opt}(i_1,j_1) \right\} \right\} \\
\Delta^\textsc{Type-2A}_\textsc{4opt}(i_2,j_2) &= \Delta^\textsc{C}_\textsc{2opt}(i_2,j_2) + \min_{i_1 \in \{0,\dots,i_2-1\}} \left\{ \min_{j_1 \in \{i_2+1,\dots,j_2-1\}} \left\{ \Delta^\textsc{D}_\textsc{2opt}(i_1,j_1) \right\} \right\} \\
\Delta^\textsc{Type-2B}_\textsc{4opt}(i_2,j_2) &= \Delta^\textsc{D}_\textsc{2opt}(i_2,j_2) + \min_{i_1 \in \{0,\dots,i_2-1\}} \left\{ \min_{j_1 \in \{i_2+1,\dots,j_2-1\}} \left\{ \Delta^\textsc{C}_\textsc{2opt}(i_1,j_1) \right\} \right\}.
\end{align}

The efficient evaluation of the best \textsc{4-opt} move for each fixed edge pair $i_2$ and $j_2$ depends on our capacity to efficiently evaluate the following quantities for $\textsc{X} \in \{\textsc{C},\textsc{D}\}$:
\begin{equation}
\Phi^\textsc{X}[i_2,j_2] =  \min_{i_1 \in \{0,\dots,i_2-1\}}  \left\{ \min_{j_1 \in \{i_2+1,\dots,j_2-1\}} \Delta^\textsc{X}_\textsc{2opt}(i_1,j_1) \right\}. \label{minmin} 
\end{equation}
By reverting the order of the minimums, observe that:
\begin{align}
\Phi^\textsc{X}[i_2,j_2] &=  \min_{j_1 \in \{i_2+1,\dots,j_2-1\}} \left\{ \min_{i_1 \in \{0,\dots,i_2-1\}}  \Delta^\textsc{X}_\textsc{2opt}(i_1,j_1) \right\} =  \min_{j_1 \in \{i_2+1,\dots,j_2-1\}} \Phi^\textsc{X}_\textsc{sub}[i_2,j_1],\\
\text{where } \Phi^\textsc{X}_\textsc{sub}[i_2,j_1] &= \min_{i_1 \in \{0,\dots,i_2-1\}}  \Delta^\textsc{X}_\textsc{2opt}(i_1,j_1).
\end{align}
Finally, observe that the $\cO(n^2)$ auxiliary variables $\Phi^\textsc{X}_\textsc{sub}$ and $\Phi^\textsc{X}$ can be calculated in $\cO(n^2)$ time with the following recursions:
\begin{equation}
\Phi^\textsc{X}_\textsc{sub}[i_2,j_1] = 
\begin{cases}
\Delta^\textsc{X}_\textsc{2opt}(0,j_1) & \text{if } i_2 = 1 \\
\min \{ \Phi^\textsc{X}_\textsc{sub}[i_2-1,j_1], \Delta^\textsc{X}_\textsc{2opt}(i_2-1,j_1)  \}  & \text{otherwise.}
\end{cases}
\label{eqaux1}
\end{equation}

\begin{equation}
\Phi^\textsc{X}[i_2,j_2] = 
\begin{cases}
\Phi^\textsc{X}_\textsc{sub}[i_2,i_2+1] & \text{if } j_2 = i_2 + 2 \\
\min \{ \Phi^\textsc{X}[i_2,j_2-1], \Phi^\textsc{X}_\textsc{sub}[i_2,j_2-1] \}  & \text{otherwise.}
\end{cases}
\label{eqaux2}
\end{equation}

Therefore, the values of these auxiliary variables can be preprocessed in $\cO(n^2)$ time.
After doing this, each evaluation of Equation~(\ref{minmin}) takes constant time, leading to a neighborhood exploration in~$\cO(n^2)$ time. The complete evaluation process is presented in Algorithm~\ref{algo:4opt}. For brevity, this pseudo-code only presents the calculation of the best-move cost, but the positions associated with the best move can be likewise recorded.

\begin{figure*}[htbp]
\centering
\begin{minipage}{0.9\textwidth}
\begin{algorithm}[H]
\IncMargin{1.5em}
\linespread{1.2}\selectfont
\caption{Efficient exploration of the \textsc{4-opt} neighborhood \label{algo:4opt}}
$\Delta_\textsc{Best} \gets 0$ \;
\For{$j_1 \in \{2,\dots,2n-1\}$}
{
\For{$i_2 \in \{1,\dots,j_1-1\}$}
{
\For{$\textsc{X} \in \{\textsc{C},\textsc{D}\}$}
{
$
\Phi^\textsc{X}_\textsc{sub}[i_2,j_1] \gets  
\begin{cases}
\Delta^\textsc{X}_\textsc{2opt}(0,j_1) & \text{if } i_2 = 1 \\
\min \{ \Phi^\textsc{X}_\textsc{sub}[i_2-1,j_1], \Delta^\textsc{X}_\textsc{2opt}(i_2-1,j_1)  \}  & \text{otherwise.}
\end{cases}
$\;
}
}
}
\For{$i_2 \in \{1,\dots,2n-2\}$}
{
\For{$j_2 \in \{i_2+2,\dots,2n\}$}
{
\For{$\textsc{X} \in \{\textsc{C},\textsc{D}\}$}
{
$
\Phi^\textsc{X}[i_2,j_2] \gets 
\begin{cases}
\Phi^\textsc{X}_\textsc{sub}[i_2,i_2+1] & \text{if } j_2 = i_2 + 2 \\
\min \{ \Phi^\textsc{X}[i_2,j_2-1], \Phi^\textsc{X}_\textsc{sub}[i_2,j_2-1] \}  & \text{otherwise.}
\end{cases}
$\;
}
$\Delta^\textsc{Type-1}_\textsc{4opt}  \ \,  (i_2,j_2)\gets  \Delta^\textsc{D}_\textsc{2opt}(i_2,j_2) + \Phi^\textsc{D}[i_2,j_2]$ \;
$\Delta^\textsc{Type-2A}_\textsc{4opt}(i_2,j_2) \gets  \Delta^\textsc{C}_\textsc{2opt}(i_2,j_2) + \Phi^\textsc{D}[i_2,j_2]$ \;
$\Delta^\textsc{Type-2B}_\textsc{4opt}(i_2,j_2) \gets  \Delta^\textsc{D}_\textsc{2opt}(i_2,j_2) + \Phi^\textsc{C}[i_2,j_2]$ \;
Check precedence constraints feasibility \;
Record the overall best feasible improving move \;
}
}
\end{algorithm}
\end{minipage}
\end{figure*}

Finally, to satisfy precedence constraints between pickup-and-delivery pairs, we include a feasibility check (Line 17 of Algorithm~\ref{algo:4opt}) which evaluates whether the best-found move of Type-1, Type-2A and Type-2B for each given $i_2$ and $j_2$ is feasible. For each move type, the conditions listed in Table~\ref{tab:conditions-feasibility} are checked in $\cO(1)$ time, using auxiliary variables which can be preprocessed in $\cO(n^2)$ time prior to move evaluations, and which state, for each sequence $\pi = \sigma_{[i,j]}$, whether this sequence is reversible (in which case $\textsc{Rev}(\pi) = \textsc{True}$) and which give the largest index $\textsc{Last}(\pi)$ of a pickup node that does not belong to $\pi$ but whose associated delivery belongs to it.

\begin{table}[htbp]
    \small
    \renewcommand{\arraystretch}{1.25}
	\caption{Feasibility conditions for \textsc{4-Opt} moves subject to pickup-and-delivery constraints \label{tab:conditions-feasibility}}
\scalebox{0.92}{
	\begin{tabular}{|c|p{10.5cm}|c|}
    \hline
    & \multicolumn{1}{c|}{\textbf{Feasibility Condition}} & \textbf{Efficient Evaluation} \\
	\hline
	\textbf{Type 1}  & $\pi_4$ cannot contain deliveries whose pickups are located in $\pi_2$ or $\pi_3$ & $\textsc{Last}(\pi_4) < |\pi_1|$  \\
            & $\pi_3$ cannot contain deliveries whose pickups are located in $\pi_2$ & $\textsc{Last}(\pi_3) < |\pi_1|$ \\
    \hline
    \textbf{Type 2A} & $\pi_3$ and $\pi_4$ cannot contain deliveries whose pickups are located in $\pi_2$ & $\max \{ \textsc{Last}(\pi_3), \textsc{Last}(\pi_4) \} < |\pi_1|$ \\
            & $\pi_3$ and $\pi_4$ should be reversible & $\textsc{Rev}(\pi_3) \land  \textsc{Rev}(\pi_4)$ \\
    \hline
    \textbf{Type 2B} & $\pi_4$ cannot contain deliveries whose pickups are located in $\pi_2$ or $\pi_3$ & $\textsc{Last}(\pi_4) < |\pi_1|$  \\
            & $\pi_2$ and $\pi_3$ should be reversible & $\textsc{Rev}(\pi_3) \land  \textsc{Rev}(\pi_4)$ \\
	\hline
	\end{tabular}
}
\end{table}%

The evaluation of the $\textsc{Rev}$ and $\textsc{Last}$ values for all sequences $\sigma_{[i,j]}$ can be done in $\cO(n^2)$ time using the following simple recursive calculation:
\begin{equation}
\textsc{Rev}(\sigma_{[i,j]}) =
\begin{cases}
\textsc{True} & \text{if } i=j \\
\textsc{Rev}(\sigma_{[i,j-1]}) \land ( \sigma(j)-n \notin \sigma_{[i,j-1]} ) & \text{otherwise},
\end{cases}
\end{equation}
\begin{equation}
\textsc{Last}(\sigma_{[i,j]}) =
\begin{cases}
\textsc{Pick}(\sigma_{[i,j]},j) & \text{if } i=j \\
\max \{\textsc{Last}(\sigma_{[i,j-1]}), \textsc{Pick}(\sigma_{[i,j]},j) \} & \text{otherwise},
\end{cases}
\end{equation}
where $\textsc{Pick}(\pi,j)$ returns the position of the pickup service associated to a delivery $\sigma(j)$ if this position does not belong to the subsequence $\pi$, and $-1$ otherwise. This function also returns $-1$ if $\sigma(j)$ is a pickup.

Overall, we achieve a complete neighborhood exploration in $\cO(n^2)$ time and space including feasibility checks. However, note that the feasibility is only evaluated for the best move of Type-1, Type-2A, and Type-2B for the pair $i_2$ and $j_2$. If this best move is infeasible, then no other move associated with this pair is applied. This assumption seems necessary to avoid a high computation effort in the feasibility checking process.

\subsection{Balas-Simonetti Neighborhood}

\cite{Balas2001} studied a TSP variant with additional constraints on the relative customer positions. This problem can be solved in polynomial time, and it leads to an efficient neighborhood search for the classical TSP. It is also possible to extend the same principle for the PDTSP, as discussed in this section.

Let $k$ be a positive integer, and let $\sigma$ be an incumbent solution. The \textsc{BS($k$)} neighborhood for the PDTSP is formed of all solutions $\pi$ satisfying the following two constraints:
\begin{align}
& \pi(i) < \pi(i+n) & \forall \, i \in \{1,\dots,n\}, \label{eq:predbs} \\
& \pi(i) < \pi(j) & \forall \, i \in V, j \in V : \sigma(i) + k \leq \sigma(j). \label{eq:bs}
\end{align}
In other words, it contains all permutations $\pi$ respecting precedence constraints between pickups and deliveries~(Equation~\ref{eq:predbs}) such that no service $i$ occurs after another service $j$ in $\pi$ if service~$j$ was planned $k$ positions later than $i$ in $\sigma$~(Equation~\ref{eq:bs}). Parameter $k$ therefore represents a degree of flexibility from the incumbent solution~$\sigma$.

In a similar way as in \cite{Balas2001}, the search for a best permutation $\pi$ satisfying Constraints~(\ref{eq:predbs}--\ref{eq:bs}) reduces to a shortest path problem in an auxiliary acyclic graph $G^* = (V^*, E^*)$. This graph is arranged into layers in such a way that only successive layers are connected by arcs. There are $2n+2$ layers in the case of the PDTSP. The first and last layers correspond to the departure and arrival at the depot $0$ and contain a single node. The other layers correspond to the $2n$ visits in $\pi$. Each such layer contains several nodes that represent combinations of anticipated and delayed visits to clients in relation to $\sigma$. More specifically, each node in $V^*$ represents a quadruplet $(i,v,S^-(\sigma,i),S^+(\sigma,i))$, where $i$ is the current position in $\pi$; $v$ is the visit (depot, pickup or delivery) allocated to position $i$ of $\pi$; $S^-(\sigma,i)$ is the set of visits that are allocated at position $i$ or later in $\sigma$ and allocated at position $i$ or before in $\pi$; and, finally, $S^+(\sigma,i)$ is the set of visits allocated before position~$i$ in~$\sigma$ that will be allocated at position~$i$ or after in~$\pi$.

Observe that sets $S^-(\sigma,i)$ and $S^+(\sigma,i)$ give enough information to restrict the search to solutions satisfying Constraints~\eqref{eq:predbs}~and~\eqref{eq:bs}.
Constraint~\eqref{eq:predbs} is ensured by restricting the search to nodes such that $S^-(\sigma,i)$ does not contain any delivery whose associated pickup belongs to \mbox{$S^+(\sigma,i) \cup v$} and~$S^+(\sigma,i)$ does not contain any pickup whose associated delivery belongs to \mbox{$S^-(\sigma,i) \cup v$}. Moreover, Constraint~(\ref{eq:bs}) is ensured by restricting the search to nodes such that all visits in $S^+(\sigma,i)$ are located no more than $k$ positions earlier than $v$ in $\sigma$.

\begin{figure}[bp]
\centering
\includegraphics[width=.9\linewidth]{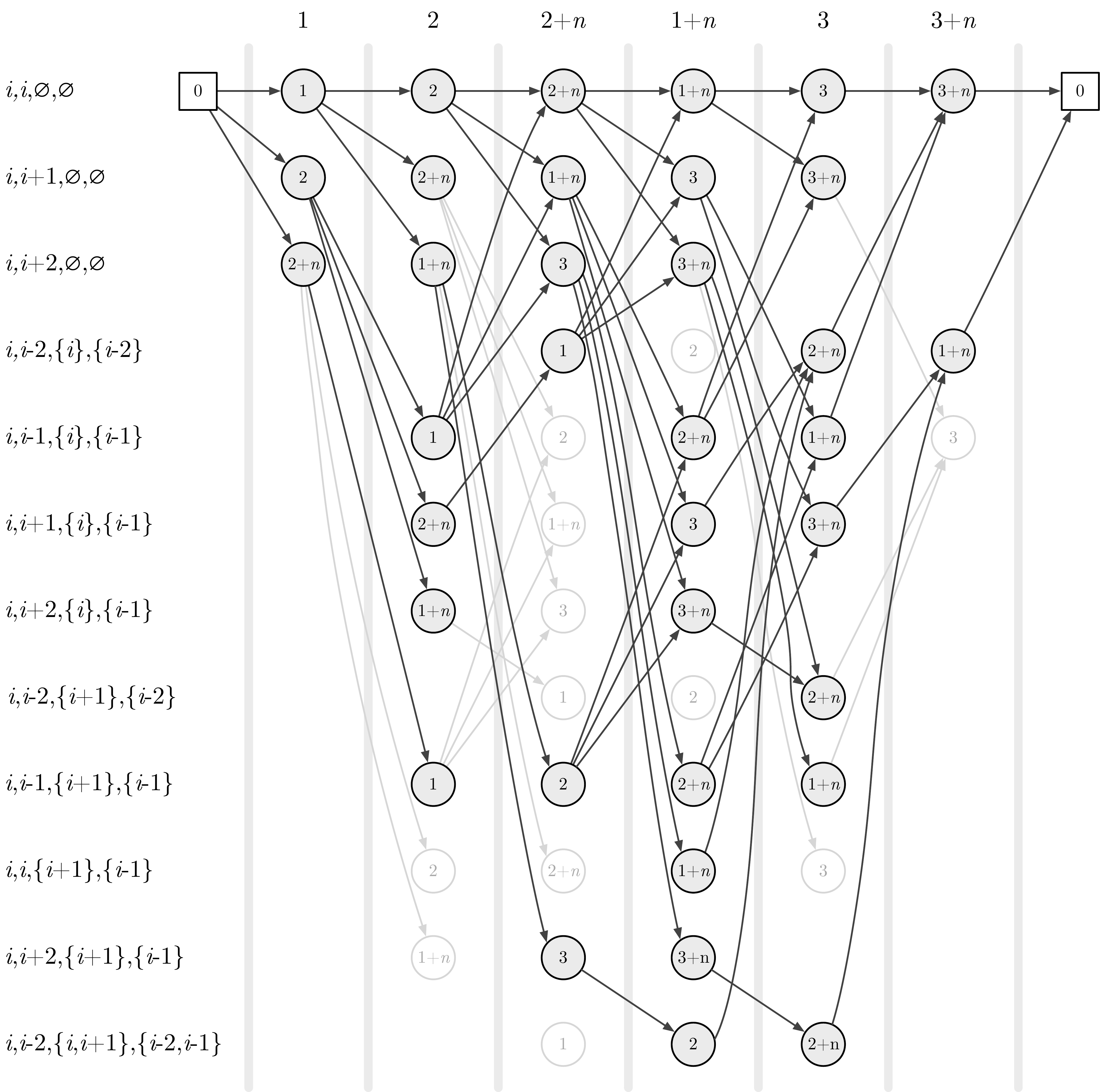}
\caption{Graph $G^*$ with $k=3$ for $\sigma = (0,1, 2, 2+n, 1+n, 3, 3+n,0)$}
\label{fig:bs}
\end{figure}

Figure \ref{fig:bs} displays $G^*$ on a simple example with three pickup-and-delivery pairs ($n=3$), starting with solution $\sigma = (0,1, 2, 2+n, 1+n, 3, 3+n,0)$ with $k=3$.
The nodes in gray represent states that are eliminated due to the pickup-and-delivery precedence constraints (Condition~\ref{eq:predbs}). For each node, we indicate the service associated with it (i.e., $v$ in the quadruplet), and on the left of the figure, we list the possible $S^-(\sigma,i)$ and $S^+(\sigma,i)$ states. The first and last layers represent the depot. The second layer, for example, contains a node corresponding to a direct visit to either $1$, $2$, or $2+n$. The third layer contains nodes representing visits to $1$, $2$, $1+n$, or $2+n$ with different combinations of anticipated and delayed visits.

Forming this graph and finding the shortest path is done in $\mathcal{O}(k^2 2^{k-2} n)$ time and $\mathcal{O}(k 2^{k-2} n)$ space. In other words, the complexity of this neighborhood exploration is linear in the number of pickup-and-delivery pairs and exponential in $k$. Notably, when $k$ is fixed, we obtain a polynomial-time algorithm to locate the best solution in this large neighborhood.

\section{Computational Experiments} 
\label{chap:5}

Most of the neighborhoods presented in this paper are either new or explored with new search strategies. The overarching goal of our experimental campaign is to measure the performance and impact of these methodological building blocks. To be representative of a broad range of methods, we conduct experiments with two structurally different heuristics: the ruin-and-recreate (RR) algorithm of \cite{Veenstra2017}, which is the state-of-the-art for the PDTSP to date, and the hybrid genetic search (HGS) of \cite{Vidal2012,Vidal2012b}, a classical population-based metaheuristic which has demonstrated good performance on a large family of node and arc routing problems \citep[see, e.g.,][]{Vidal2014,Vidal2017b}.
Notably, the principle that we proposed for an efficient exploration of the \textsc{Relocate Pair} neighborhood (in Section~\ref{relocate-pair}) can also be applied to efficiently evaluate candidate insertions in the context of a large neighborhood search such as RR. As this permits considerable reductions of complexity and computational effort or, alternatively, an increased number of iterations evaluated in a given time, we will also present comparative analyses demonstrating the impact of this improvement in the RR methodology.

Our experiments are conducted on three main groups of instances: (i) classic PDTSP instances with 11 to 493 visits, (ii) new challenging instances derived from the X set of \citet{Uchoa2017} with 101 to 1001 visits, and (iii) the ``Grubhub'' instance set (available at \url{https://github.com/grubhub/tsppdlib}) connected to meal-delivery applications. The latter instances only include up to 31 visits, but nonetheless remain challenging for exact methods. Due to their repeated solution in a delivery-on-demand setting, it is essential to consistently solve them to near-optimality within the shortest possible time, ideally a few milliseconds.

In Section \ref{sec:exp:meta}, we describe the local search and metaheuristic frameworks used in our experiments. In Section \ref{sec:exp:setup}, we detail the experimental setup, benchmark instances and parameter-calibration methodology. In Section~\ref{sec:exp:comparisons}, we compare the performance of the considered metaheuristics as well as the impact of the fast evaluation of insertions in the RR. In Section~\ref{sec:exp:sensitivity}, we measure the sensitivity of HGS to changes in the choices of neighborhoods. Finally, in Section \ref{sec:exp:mealdelivery} we report results on the Grubhub data sets, which require very short solution times.

\subsection{Search Strategies}
\label{sec:exp:meta}

We first describe the local search built on the proposed neighborhoods, and then proceed with a description of the RR and HGS metaheuristics. Unless specified, these metaheuristics were used with their original parameters listed in~\cite{Vidal2012} and~\cite{Veenstra2017}.

\subsubsection{Local Search.}
\label{exp:meta:1}
The proposed local search algorithm is built on the six neighborhoods presented in Section~\ref{sec:neighborhoods}.
As illustrated in Algorithm~\ref{Algo:LS}, the exploration of the neighborhoods is divided into two alternating phases.
In a first phase, it explores the \textsc{Relocate Pair}, \textsc{2-Opt} and \textsc{Or-Opt} neighborhoods, dividing the evaluation of the moves among the different pickup-and-delivery pairs $(x,x+n)$ and applying the best move for each such pair. Next, it performs a second-phase search that consists of iteratively finding and applying the best \textsc{2k-Opt}, \textsc{4-Opt} or \textsc{Balas-Simonetti} move. The search stops whenever all those neighborhoods have been explored without improvement.

\begin{figure*}[htbp]
\centering
\begin{minipage}{\textwidth}
\begin{algorithm}[H]
\IncMargin{1.5em}
\linespread{1.2}\selectfont
\caption{General structure of the local search}
\label{Algo:LS}
$S \gets \text{Initial Solution}$ \;
\While{\emph{Improvements Found}}
{
\For{\emph{all pickup-and-delivery requests $(x,x+n)$ in random order}\hfill \tcp{$n$ iterations}}
{
Find and apply, if improving, the best move on $S$ among:\nonl

\hspace*{1em} $\bullet$ all \textsc{Relocate Pair} moves for $(x,x+n)$; \hfill \tcp{$\cO(n)$ time} \nonl

\hspace*{1em} $\bullet$ all \textsc{2-Opt} moves starting at the position of $x$ or $x+n$; \hfill \tcp{$\cO(n)$ time} \nonl

\hspace*{1em } $\bullet$ all \textsc{Or-Opt} moves starting at the position of $x$ or $x+n$. \hfill \tcp{$\cO(k_\textsc{or}n)$ time}
}
Find and apply, if improving, the best move on $S$ in \textsc{2k-Opt};
\textsc{4-Opt};\label{lineLarge}\nonl

and \textsc{Balas-Simonetti} neighborhoods. \hfill \tcp{$\cO(k^2 2^{k-2} n + n^2)$ time}
}
Return~$S$ \;
\end{algorithm}
\end{minipage}
\end{figure*}

\subsubsection{Ruin-and-Recreate.} 
\label{exp:meta:2}
The RR algorithm of \cite{Veenstra2017} is based on a simple iterative search process. It starts from an initial solution built with the basic greedy insertion procedure of \cite{Ropke2006b}, which consists of shuffling the order of the pickup-and-delivery pairs and inserting each of them in their best position.
Then, it iteratively destroys and reconstructs the incumbent solution. Four destruction operators are used.
\begin{itemize}[leftmargin=*]
\item \textsc{Random removal:} Remove $q$ random pickup-and-delivery services.
\item \textsc{Worst removal:} Calculate the gain associated with the removal of each pickup-and-delivery service. List the services by decreasing gain. Remove $q$ random services, giving a higher probability to elements that are early in the list.
\item \textsc{Block removal:} Iteratively select a random pickup-and-delivery service $(x,x+n)$. Remove all the services containing a pickup or a delivery in the sequence starting at $x$ and finishing at $x+n$. Stop as soon as $q$ services have been removed.
\end{itemize}
The number of pickup-and-delivery services $q$ selected for removal is each time randomly selected, with uniform probability, between $\min\{30,0.20n\}$ and $\min\{50,0.55n\}$. Reconstruction is done with the same greedy insertion heuristic as in the initialization phase, considering only the removed services. The resulting solution is accepted as a new incumbent solution in case of improvement, and otherwise subject to a probabilistic acceptance using Metropolis criterion (as in the simulated annealing algorithm).

In \cite{Veenstra2017}, each evaluation of the best insertion position for a pickup-and-delivery pair is done by enumeration (on the pickup and delivery insertion position) and therefore takes $\cO(n^2)$ time, leading to a total complexity of $\cO(q n^2)$ for each reconstruction phase. Notably, applying the same strategy as in the \textsc{Relocate Pair} through Equations~(\ref{eq:rp:1}--\ref{eq:rp:2}) permits to cut down computational complexity to $\cO(q n)$ time for the complete reconstruction phase. This significantly speeds up the overall RR approach since the reconstruction operator is the main bottleneck of the method. We will show the impact of this improvement in Section~\ref{sec:exp:comparisons}.

\subsubsection{Hybrid Genetic Search.}
\label{exp:meta:3}
Our HGS for the PDTSP follows the original search strategy of \cite{Vidal2012} and \cite{Vidal2022} with some adaptations and simplifications. Since a solution of the PDTSP consists of a single tour, the algorithm maintains a single population of feasible solutions represented as permutations of the pickup and delivery vertices.
To produce each new solution, our HGS selects two parents by binary tournament, applies a crossover and a mutation operator to generate a new child, which is then improved using the proposed local search (see Section~\ref{exp:meta:1}) and included in the~population. To reduce the computational effort, the complete local search with the large neighborhoods is only used for a random subset of the generated individuals selected with probability~$p_\textsc{large}$. For the other individuals, the local search relies exclusively on \textsc{Relocate Pair}, \textsc{2-Opt} and \textsc{Or-Opt} (i.e., Algorithm~\ref{Algo:LS} without Line~\ref{lineLarge}).

\begin{figure}[htbp]
\centering
\includegraphics[width=0.9\linewidth]{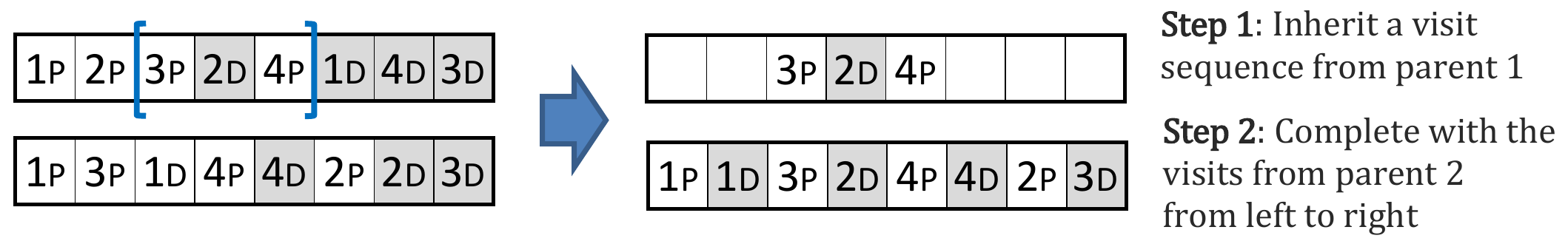}
\caption{Illustration of the \textsc{LOX} crossover}
\label{fig:LOX}
\end{figure}

We use the linear order crossover (LOX -- \citealt{Falkenauer1991}) illustrated in Figure~\ref{fig:LOX}. This operator consists of copying a random sequence of consecutive visits from the first parent into the child and completing the missing visits from the left to the right based on the visit order of the second parent. The crossover is followed by a mutation which consists of finding and applying the best \textsc{4-Opt} move of Type~1, regardless of its feasibility. These two operations may lead to precedence-constraints violations. Therefore, in a last step, HGS applies the best \textsc{Relocate Pair} move for each pickup-and-delivery pair violating precedence constraints to obtain the best feasible reinsertion positions and achieve feasibility. This mutation operator successfully introduces promising edges and contributes to the diversification of the search.

As usual in HGS, we evaluate the \emph{biased fitness} $BF(\sigma)$ of each solution as a linear combination of its rank in terms of cost and its rank in terms of contribution to the population diversity.
\begin{equation}
BF(\sigma) = R_\textsc{c}(\sigma) + \left( 1 - \frac{\mu^\textsc{Elite}}{|\mathcal{P}|} \right) \times R_\textsc{d}(\sigma) \label{eq-bf}
\end{equation}
This fitness measure is used for parents and survivors selections. 
The parameter $\mu^\textsc{Elite}$ represents a number of elite individuals which are guaranteed to survive. It has been set to $\mu^\textsc{Elite} = 1$ in our experiments to promote a diversified search and preserve the best solution. $R_\textsc{c}(\sigma)$ represents the rank of solution $\sigma$ within the population $\mathcal{P}$ in terms of objective value, and $R_\textsc{d}(\sigma)$ represents its rank in terms of diversity contribution. The diversity contribution of each solution is computed as its average distance to the two closest other solutions in $\mathcal{P}$. To calculate the distance $\delta(\sigma,\sigma')$ between two solutions~$\sigma$ and~$\sigma'$, we use the Jaccard distance \citep{Levandowsky1971} between the sets $E(\sigma)$ and~$E(\sigma')$ containing the pairs of consecutive visit locations within each solution (i.e., considering two visits in the same location as identical):
\begin{equation}
\delta(\sigma,\sigma') = \frac{|E(\sigma) \cup E(\sigma')| - |E(\sigma) \cap E(\sigma')| }{|E(\sigma) \cup E(\sigma')|}. \label{eq:jaccard}
\end{equation}
This choice of distance metric is due to the fact that PDTSP applications and datasets commonly include  consecutive pickups or deliveries from different services at the same locations. It avoids distance over-estimations for pairs of symmetrical solutions with different visit permutations in the same location.

Finally, the population is managed to contain between $\mu$ and $\mu +\lambda$ solutions. Initially, it contains random solutions generated by the same randomized greedy algorithm as in Section~\ref{exp:meta:2}. Each time the maximum population size is attained, the method iteratively removes the worst solution in terms of biased fitness until reaching the initial population size of $\mu$. 

\subsection{Experimental Setup}
\label{sec:exp:setup}

All our computational experiments were conducted on a single thread of an AMD Rome 7532 2.40GHz CPU running Linux CentOS 8.4. The algorithms were implemented in \textsc{C++}, and compiled with \textsc{G++} v7.4.0 using option \textsc{-O3}.

The classical benchmark instances of the PDTSP are classified into three sets. The first two classes of instances (RBO) were introduced in \cite{Renaud2000a} and the third class (DRCL) was introduced in \cite{Dumitrescu2010}. The instances of the second class are very peculiar since, by design, the optimal PDTSP solution on each instance matches that of the TSP, while the instances of the third class are fairly small. For these reasons, we focus our studies on the instances of the first class, also including computational results for the other classes in the electronic companion.
The first class includes 108 instances containing 51 to 493 visits. These instances are derived from 36 TSPLIB data sets \citep{Reinelt1991} and subdivided into three groups called~A, B, and C. For each TSP data set, three PDTSP instances were generated by randomly picking a delivery point for each (pickup) vertex. In groups A and B, one random delivery point was selected for each pickup among the five and ten closest vertices, respectively. In group C, there is no such limitation, and therefore the delivery vertex was freely selected among all other vertices. To obtain a wider diversity of challenging data sets, we also derived 300 new instances for the PDTSP from the X instance set of \cite{Uchoa2017}. To do so, we retained the coordinates of the depot and customers for each of the 100 original instances, and we created pickup and delivery pairs using the same process as previously described, thereby creating three instances (A, B and C). We use the same rounding conventions as in the set X for distance calculations, i.e., we round the distance values to the nearest integer.\\

\noindent
\textbf{Calibration of the neighborhood search.}
Our neighborhood search includes three main parameters: the probability of using the large neighborhoods ($p_\textsc{large}$), the maximum size of the \textsc{Or-Opt} sequences ($k_\textsc{or}$), and the \textsc{Balas-Simonetti} parameter ($k_\textsc{bs}$).
To calibrate these parameters, a full factorial experiment would take excessive computational effort. Therefore, we conducted preliminary analyses to identify a suitable level for $p_\textsc{large}$, leading us to a value of \mbox{$p_\textsc{large} = 0.1$}. Next, we calibrated $k_\textsc{or}$ and $k_\textsc{bs}$, and lastly we reevaluated our choice of $p_\textsc{large}$ to make sure that it remained meaningful. For this calibration, we used the 300 new instances of set X and ran HGS up until a fixed time limit of $T_\textsc{max} = N$ seconds, where $N$ represents the total number of visits in the instance. We performed ten runs with different random seeds for each instance and algorithm configuration. Table~\ref{tab:lscalibration} reports the performance of each configuration, measured as an average error gap over all instances and runs. For each instance, the gap is calculated relatively to the best-known solution (BKS) collected from all runs as $\text{Gap}(\%) = 100 \times (z_\textsc{Sol} - z_\textsc{bks}) / z_\textsc{bks}$. The table also includes a color scale to highlight better configurations, and the overall best configuration is marked in~boldface.

\begin{table}[htbp]
	\caption{Calibration of the local search parameters}	
	\label{tab:lscalibration}
    \setlength\tabcolsep{12pt}
    \renewcommand{\arraystretch}{1.2}
    \centering
	\scalebox{0.78}
	{
	\begin{tabular}{|l|c|c|c|c|c|c|}
		\hline
		&  $k_\textsc{bs} = 1$  &
		   $k_\textsc{bs} = 2$  & 
		   $k_\textsc{bs} = 3$  & 
		   $k_\textsc{bs} = 4$  & 
		   $k_\textsc{bs} = 5$  & 
		   $k_\textsc{bs} = 6$  \\
		\hline
        $k_\textsc{or} = 1$ & \cellcolor[rgb]{0.85581, 0.888601, 0.097452}0.675 & \cellcolor[rgb]{0.82494, 0.88472, 0.106217}0.667 & \cellcolor[rgb]{0.886271, 0.892374, 0.095374}0.680 & \cellcolor[rgb]{0.935904, 0.89857, 0.108131}0.688 & \cellcolor[rgb]{0.906311, 0.894855, 0.098125}0.686 & \cellcolor[rgb]{0.964894, 0.902323, 0.123941}0.696 \\ 
        $k_\textsc{or} = 2$ & \cellcolor[rgb]{0.699415, 0.867117, 0.175971}0.590 & \cellcolor[rgb]{0.762373, 0.876424, 0.137064}0.592 & \cellcolor[rgb]{0.730889, 0.871916, 0.156029}0.591 & \cellcolor[rgb]{0.79376, 0.880678, 0.120005}0.598 & \cellcolor[rgb]{0.668054, 0.861999, 0.196293}0.588 & \cellcolor[rgb]{0.636902, 0.856542, 0.21662}0.586 \\ 
        $k_\textsc{or} = 3$ & \cellcolor[rgb]{0.496615, 0.826376, 0.306377}0.523 & \cellcolor[rgb]{0.555484, 0.840254, 0.269281}0.531 & \cellcolor[rgb]{0.468053, 0.818921, 0.323998}0.507 & \cellcolor[rgb]{0.575563, 0.844566, 0.256415}0.540 & \cellcolor[rgb]{0.525776, 0.833491, 0.288127}0.529 & \cellcolor[rgb]{0.606045, 0.850733, 0.236712}0.544 \\ 
        $k_\textsc{or} = 5$ & \cellcolor[rgb]{0.360741, 0.785964, 0.387814}0.489 & \cellcolor[rgb]{0.335885, 0.777018, 0.402049}0.486 & \cellcolor[rgb]{0.311925, 0.767822, 0.415586}0.479 & \cellcolor[rgb]{0.386433, 0.794644, 0.372886}0.486 & \cellcolor[rgb]{0.412913, 0.803041, 0.357269}0.491 & \cellcolor[rgb]{0.440137, 0.811138, 0.340967}0.501 \\ 
        $k_\textsc{or} = 7$ & \cellcolor[rgb]{0.232815, 0.732247, 0.459277}0.453 & \cellcolor[rgb]{0.196571, 0.711827, 0.479221}0.448 & \cellcolor[rgb]{0.214, 0.722114, 0.469588}0.451 & \cellcolor[rgb]{0.252899, 0.742211, 0.448284}0.455 & \cellcolor[rgb]{0.288921, 0.758394, 0.428426}0.458 & \cellcolor[rgb]{0.274149, 0.751988, 0.436601}0.457 \\ 
        $k_\textsc{or} = 10$ & \cellcolor[rgb]{0.143303, 0.669459, 0.511215}0.428 & \cellcolor[rgb]{0.180653, 0.701402, 0.488189}0.439 & \cellcolor[rgb]{0.153894, 0.680203, 0.504172}0.437 & \cellcolor[rgb]{0.12478, 0.640461, 0.527068}0.423 & \cellcolor[rgb]{0.166383, 0.690856, 0.496502}0.438 & \cellcolor[rgb]{0.134692, 0.658636, 0.517649}0.427 \\ 
        $k_\textsc{or} = 15$ & \cellcolor[rgb]{0.136408, 0.541173, 0.554483}0.416 & \cellcolor[rgb]{0.174274, 0.445044, 0.557792}0.411 & \cellcolor[rgb]{0.147607, 0.511733, 0.557049}0.415 & \cellcolor[rgb]{0.12138, 0.629492, 0.531973}0.422 & \cellcolor[rgb]{0.151918, 0.500685, 0.557587}0.415 & \cellcolor[rgb]{0.119699, 0.61849, 0.536347}0.421 \\ 
        $k_\textsc{or} = 20$ & \cellcolor[rgb]{0.140536, 0.530132, 0.555659}0.416 & \cellcolor[rgb]{0.165117, 0.467423, 0.558141}0.412 & \cellcolor[rgb]{0.208623, 0.367752, 0.552675}0.408 & \cellcolor[rgb]{0.128087, 0.647749, 0.523491}0.424 & \cellcolor[rgb]{0.225863, 0.330805, 0.547314}0.396 & \cellcolor[rgb]{0.183898, 0.422383, 0.556944}0.409 \\ 
        $k_\textsc{or} = 30$ & \cellcolor[rgb]{0.128729, 0.563265, 0.551229}0.418 & \cellcolor[rgb]{0.160665, 0.47854, 0.558115}0.415 & \cellcolor[rgb]{0.231674, 0.318106, 0.544834}\textbf{0.395} & \cellcolor[rgb]{0.197636, 0.391528, 0.554969}0.408 & \cellcolor[rgb]{0.192357, 0.403199, 0.555836}0.409 & \cellcolor[rgb]{0.132444, 0.552216, 0.553018}0.418 \\ 
        $k_\textsc{or} = 40$ & \cellcolor[rgb]{0.214298, 0.355619, 0.551184}0.407 & \cellcolor[rgb]{0.120565, 0.596422, 0.543611}0.420 & \cellcolor[rgb]{0.179019, 0.433756, 0.55743}0.410 & \cellcolor[rgb]{0.187231, 0.414746, 0.556547}0.409 & \cellcolor[rgb]{0.203063, 0.379716, 0.553925}0.408 & \cellcolor[rgb]{0.119512, 0.607464, 0.540218}0.420 \\ 
        $k_\textsc{or} = 50$ & \cellcolor[rgb]{0.220057, 0.343307, 0.549413}0.402 & \cellcolor[rgb]{0.125394, 0.574318, 0.549086}0.420 & \cellcolor[rgb]{0.169646, 0.456262, 0.55803}0.411 & \cellcolor[rgb]{0.143343, 0.522773, 0.556295}0.415 & \cellcolor[rgb]{0.122606, 0.585371, 0.546557}0.420 & \cellcolor[rgb]{0.15627, 0.489624, 0.557936}0.415 \\
		\hline
	\end{tabular}%
	}
\end{table}%

As visible in these experiments, parameter $k_\textsc{or}$ has a noticeable impact on the performance of the overall method. This parameter limits the size of the \textsc{Or-Opt} neighborhood, which is especially important for the PDP as it permits structural changes to the solutions without reverting visit sequences. However, the exploration of this neighborhood is a time bottleneck when $k_\textsc{or}$ is large. Our experiments show that suitable values for $k_\textsc{or}$ are located in the $\{20,\dots,50\}$ range. In contrast, parameter $k_\textsc{bs}$ has a more limited impact, though excessively small or large values (e.g., $1$ or $6$) also lead to a performance deterioration. Based on these results, we adopted as \emph{reference configuration} the values $(k_\textsc{or},k_\textsc{bs}) = (30,3)$, which achieved the best performance. As displayed in Table~\ref{tab:lscalibrationwx}, we also ran pairwise Wilcoxon tests to compare the results of the reference configuration with those of the 65 other configurations. At a significance level of $0.05/65$ (using Bonferroni correction to account for multiple testing), the reference configuration is significantly better than any other configuration with $k_\textsc{or} \leq 10$. Finally, as seen on Table~\ref{tab:lscalibrationplarge}, re-calibrating $p_\textsc{Large}$ based on the reference values of $(k_\textsc{or},k_\textsc{bs})$ confirmed our initial choice of \mbox{$p_\textsc{Large} = 0.1$}.

\begin{table}[htbp]
	\centering
	\caption{p-values of paired-samples Wilcoxon tests relative to the reference configuration}
	\label{tab:lscalibrationwx}
    \setlength\tabcolsep{8pt}
    \renewcommand{\arraystretch}{1.2}
	\scalebox{0.78}
	{
	\begin{tabular}{|l|c|c|c|c|c|c|}
		\hline
		&  $k_\textsc{bs} = 1$  &
		   $k_\textsc{bs} = 2$  & 
		   $k_\textsc{bs} = 3$  & 
		   $k_\textsc{bs} = 4$  & 
		   $k_\textsc{bs} = 5$  & 
		   $k_\textsc{bs} = 6$  \\
		\hline
        $k_\textsc{or} = 1$ & {\num{1.43335471867304E-153}} & {\num{8.48026672463656E-149}} & {\num{1.00E-155}} & {\num{9.68820717478757E-162}} & {\num{4.60197966521953E-156}} & {\num{8.87933095157177E-153}} \\ 
        $k_\textsc{or} = 2$ & {\num{2.43953330987366E-93}} & {\num{2.28637404140016E-108}} & {\num{3.09683839117272E-92}} & {\num{2.20191629671723E-101}} & {\num{1.88370281741205E-88}} & {\num{7.82064755080316E-92}} \\ 
        $k_\textsc{or} = 3$ & {\num{4.4808258137422E-50}} & {\num{3.61488190063454E-59}} & {\num{2.61127732251831E-43}} & {\num{2.19438793580457E-63}} & {\num{1.91558054114948E-55}} & {\num{1.61843464579875E-66}} \\ 
        $k_\textsc{or} = 5$ & {\num{2.38504467489973E-30}} & {\num{1.70871066468507E-28}} & {\num{7.74955076731098E-27}} & {\num{1.92040271589232E-31}} & {\num{1.73108871689606E-34}} & {\num{2.71677297952884E-38}} \\ 
        $k_\textsc{or} = 7$ & {\num{7.0776790688508E-16}} & {\num{2.15363545403289E-15}} & {\num{2.18902451038108E-14}} & {\num{7.59648631944355E-14}} & {\num{3.25210532788005E-15}} & {\num{7.46824121895275E-17}} \\ 
        $k_\textsc{or} = 10$ & {\num{1.30860973164162E-05}} & {\num{2.39119824858872E-08}} & {\num{5.01706971548481E-09}} & {\num{1.31332700262486E-05}} & {\num{1.35789722835753E-09}} & {\num{1.39768272445319E-05}} \\ 
        $k_\textsc{or} = 15$ & {\num{0.001848641831008}} & {\num{0.007708308493519}} & {\num{0.000732485584308}} & {\num{5.33653830022716E-05}} & {\num{0.00350426640439}} & {\num{0.000377267324838}} \\ 
        $k_\textsc{or} = 20$ & {\num{0.014217206572145}} & {\num{0.028305888841541}} & \num{0.168544271862022} & {\num{8.13241472537867E-05}} & \num{0.911589043333628} & \num{0.056003266112891} \\ 
        $k_\textsc{or} = 30$ & {\num{0.009519081528741}} & \num{0.058546754000831} & --- & \num{0.07938671860578} & \num{0.144772522240583} & {\num{0.008234812465381}} \\ 
        $k_\textsc{or} = 40$ & \num{0.2875164434639} & {\num{0.000711186829175}} & \num{0.068893801899332} & {\num{0.02909623914125}} & {\num{0.025646208239451}} & {\num{0.012437384908729}} \\ 
        $k_\textsc{or} = 50$ & {\num{0.047435335263802}} & {\num{0.002875476105923}} & \num{0.070814227600056} & {\num{0.010430856172022}} & {\num{0.00076558014802}} & {\num{0.000375788439338}} \\ 
		\hline
	\end{tabular}%
	}
\end{table}%

\begin{table}[htbp]
	\centering
	\caption{Final calibration of the $p_\text{Large}$ parameter}
	\label{tab:lscalibrationplarge}
    \renewcommand{\arraystretch}{1.2}
    \setlength\tabcolsep{10pt}
	\scalebox{0.8}
	{
	\begin{tabular}{|l|c|c|c|c|}
		\hline
		$p_\textsc{large}$ &  $k_\textsc{or}$  & $k_\textsc{bs}$ &  Gap(\%) & p-value \\
		\hline
        0.00 & 30 & 3 & 0.423 & \num{0.003229077825285} \\ 
        0.05 & 30 & 3 & 0.428 & \num{0.000538189477981} \\ 
        0.10 & 30 & 3 & \textbf{0.395} & --- \\ 
        0.20 & 30 & 3 & 0.422 & \num{0.002982932714182} \\ 
        0.30 & 30 & 3 & 0.427 & \num{0.000197686792092} \\ 
        0.50 & 30 & 3 & 0.434 & \num{3.02384225952432E-06} \\ 
        1.00 & 30 & 3 & 0.490 & \num{1.89360889137343E-24} \\ 
		\hline
	\end{tabular}%
		}
\end{table}%

\subsection{Performance Comparisons}
\label{sec:exp:comparisons}

The goal of this section is to establish a comparative analysis of the different PDTSP solution methods. We use new seed values to avoid any bias due to possible calibration overfit. As we re-implemented the RR approach, we can establish a comparison with the same computational environment and termination criterion. We therefore use a common time limit of $T_\textsc{max} = N$ seconds for each  instance. The cooling schedule of the RR method has also been calibrated to follow an exponential decay from the initial temperature to the final temperature within the allotted time.
Finally, we confront the new experimental results with those reported from previous studies, which serve as a baseline (i.e., to compare our re-implementation of RR with the original results and previous approaches). Overall, we report the results of the following heuristics:
\begin{itemize}
    \item RBO-Best: Best results from \cite{Renaud2000a} and \cite{Renaud2002a}, directly provided to us by the authors. We note that these solution values were registered from multiple runs (136 and 168, respectively) using different parameter values.
    \item RR-VRVC: Results of the RR algorithm of \cite{Veenstra2017} over ten runs, as reported by the authors. Experiments conducted on an Intel Core i3-2120 @3.3GHz CPU.
    \item RR and RR-fast: Results of our re-implementation of the RR algorithm over ten runs, without or with the speed-up strategy for the reconstruction step, respectively, using $T_\textsc{max} = N$.
    \item HGS: Results of the hybrid genetic search over ten runs, using $T_\textsc{max} = N$.
\end{itemize}
The CPU used in \cite{Veenstra2017} is estimated to be $1.41 \times$ slower than ours according to Passmark single-thread benchmark. Therefore, we divide the CPU time values of RR-VRVC by this coefficient to draw a reliable comparison.

Tables~\ref{tab:summaryres1}~and~\ref{tab:summaryresX} compare the performance of all of these methods on the classical instances and the new ones derived from set X, respectively.
These tables list, when available, the best solutions from RBO as well as the average percentage gap (Avg) relative to the BKS over ten runs for each instance, the best percentage gap (Best), and the average (scaled) computational time in seconds (T(s)) for each other method. The gap is calculated as $\text{Gap}(\%) = 100 \times (z_\textsc{Sol} - z_\textsc{bks}) / z_\textsc{bks}$ for each instance and solution. In Table~\ref{tab:summaryres1}, the BKS have been collected from all previous works by \citet{Renaud2000a,Renaud2002a}, \citet{Dumitrescu2010} and \citet{Veenstra2017}, therefore negative values represent improvements over the previous BKSs. In Table~\ref{tab:summaryresX}, as no previous results are available, the BKSs correspond to the best solutions ever found by any method on any run, including the calibration experiments.
The results of these tables are aggregated per range of instance size and instance type, and the best average solution quality is indicated on each line. The detailed solution values for all instances are provided in the Electronic Companion. Moreover, the source code and scripts needed to run these experiments, the new instances, as well as the solutions for each instance are openly accessible at \url{https://github.com/vidalt/PDTSP}.

\begin{table}[htbp]
	\centering
	\caption{Performance comparison on the instances of Class 1} \label{tab:summaryres1}
	\setlength\tabcolsep{6.5pt}
	\renewcommand{\arraystretch}{1.05}
	\scalebox{0.7}
	{
		\begin{tabular}{c@{\hspace*{0.2cm}}c@{\hspace*{0.2cm}}cccccccccccccccccccc}
			\toprule
			& & & RBO & & \multicolumn{3}{c}{RR-VRVC} & & \multicolumn{3}{c}{RR} & & \multicolumn{3}{c}{RR-fast}  & & \multicolumn{3}{c}{HGS} \\
			\cmidrule{4-4} \cmidrule{6-8} \cmidrule{10-12} \cmidrule{14-16} \cmidrule{18-20}
			$N$ & Type & & Best & & Avg & Best & T(s) & & Avg & Best & T(s) & & Avg & Best & T(s) & & Avg & Best & T(s) \\
			\midrule
			
			 \multirow{3}[0]{*}{50-99}
			 & A &  & 0.015 &  & 0.222 & 0.000 & 17.38 &  & 0.005 & -0.020 & 87.55 &  & -0.008 & -0.020 & 87.55 &  & \textbf{-0.023} & -0.023 & 87.55 \\ 
             & B &  & 0.000 &  & 0.128 & 0.027 & 18.54 &  & 0.022 & 0.011 & 87.55 &  & 0.009 & -0.014 & 87.55 &  & \textbf{-0.014} & -0.014 & 87.55 \\ 
             & C &  & 0.000 &  & 0.060 & 0.007 & 19.73 &  & 0.007 & 0.007 & 87.55 &  & 0.007 & 0.007 & 87.55 &  & \textbf{0.007} & 0.007 & 87.55 \\
             & & & & & & & & & & & & & & & & & & & \\
             
             \multirow{3}[0]{*}{100-199}
             & A &  & 0.401 &  & 0.353 & 0.016 & 83.59 &  & 0.000 & -0.132 & 146.00 &  & -0.120 & -0.169 & 146.00 &  & \textbf{-0.172} & -0.172 & 146.00 \\ 
             & B &  & 0.488 &  & 0.260 & 0.002 & 87.05 &  & 0.008 & -0.064 & 146.00 &  & -0.056 & -0.096 & 146.00 &  & \textbf{-0.098} & -0.098 & 146.00 \\ 
             & C &  & 0.482 &  & 0.366 & 0.000 & 93.67 &  & 0.168 & 0.012 & 146.00 &  & 0.027 & -0.061 & 146.00 &  & \textbf{-0.074} & -0.074 & 146.00 \\
             & & & & & & & & & & & & & & & & & & & \\
             
             \multirow{3}[0]{*}{200-299}
             & A &  & 1.296 &  & 0.824 & 0.104 & 364.28 &  & 0.508 & 0.101 & 254.60 &  & 0.014 & -0.413 & 254.60 &  & \textbf{-0.580} & -0.580 & 254.60 \\ 
             & B &  & 1.118 &  & 0.643 & 0.000 & 381.73 &  & 0.070 & -0.481 & 254.60 &  & -0.396 & -0.561 & 254.60 &  & \textbf{-0.605} & -0.625 & 254.60 \\ 
             & C &  & 1.618 &  & 0.691 & 0.000 & 414.33 &  & 0.099 & -0.179 & 254.60 &  & -0.229 & -0.460 & 254.60 &  & \textbf{-0.516} & -0.580 & 254.60 \\
             & & & & & & & & & & & & & & & & & & & \\
             
             \multirow{3}[0]{*}{300-499}
             & A &  & 2.974 &  & 0.852 & 0.000 & 1015.46 &  & 0.438 & -0.247 & 417.67 &  & -0.391 & -0.613 & 417.67 &  & \textbf{-0.844} & -1.030 & 417.67 \\ 
             & B &  & 2.778 &  & 0.863 & 0.000 & 1027.61 &  & 0.543 & -0.043 & 417.67 &  & -0.099 & -0.487 & 417.67 &  & \textbf{-0.847} & -1.070 & 417.67 \\ 
             & C &  & 1.698 &  & 1.159 & 0.025 & 1126.87 &  & 0.940 & -0.331 & 417.67 &  & -0.189 & -0.635 & 417.67 &  & \textbf{-1.249} & -1.484 & 417.67 \\
             & & & & & & & & & & & & & & & & & & & \\
             
             \multirow{3}[0]{*}{All} 
             & A &  & 0.836 &  & 0.462 & 0.021 & 257.66 &  & 0.145 & -0.084 & 188.50 &  & -0.112 & -0.231 & 188.50 &  & \textbf{-0.295} & -0.326 & 188.50 \\ 
             & B &  & 0.808 &  & 0.373 & 0.009 & 263.81 &  & 0.110 & -0.096 & 188.50 &  & -0.090 & -0.201 & 188.50 &  & \textbf{-0.268} & -0.307 & 188.50 \\ 
             & C &  & 0.695 &  & 0.450 & 0.006 & 287.81 &  & 0.238 & -0.073 & 188.50 &  & -0.051 & -0.192 & 188.50 &  & \textbf{-0.307} & -0.355 & 188.50 \\
             
			\midrule
			\multicolumn{2}{c}{Average} & & 0.780 &  & 0.428 & 0.012 & 269.76 &  & 0.164 & -0.084 & 188.50 &  & -0.085 & -0.208 & 188.50 &  & \textbf{-0.290} & -0.329 & 188.50 \\ 
			\bottomrule
		\end{tabular}
	}
\end{table}

\begin{table}[htbp]
	\centering
	\caption{Performance comparison on the instances of Class X} \label{tab:summaryresX}
	\setlength\tabcolsep{6.5pt}
	\renewcommand{\arraystretch}{1.05}
	\scalebox{0.74}
	{
		\begin{tabular}{c@{\hspace*{0.2cm}}c@{\hspace*{0.2cm}}ccccccccccccc}
			\toprule
			& & & \multicolumn{3}{c}{RR} & & \multicolumn{3}{c}{RR-fast}  & & \multicolumn{3}{c}{HGS} \\
			\cmidrule{4-6} \cmidrule{8-10} \cmidrule{12-14}
			$N$ & Type & & Avg & Best & T(s) & & Avg & Best & T(s) & & Avg & Best & T(s) \\
			\midrule
			
            \multirow{3}[0]{*}{100-199}
            & A &  & 0.251 & 0.132 & 147.57 &  & 0.166 & 0.072 & 147.58 &  & \textbf{0.000} & 0.000 & 147.57 \\ 
            & B &  & 0.167 & 0.044 & 147.57 &  & 0.063 & 0.010 & 147.58 &  & \textbf{0.001} & 0.000 & 147.57 \\ 
            & C &  & 0.276 & 0.156 & 147.57 &  & 0.146 & 0.093 & 147.57 &  & \textbf{0.013} & 0.003 & 147.57 \\ 
            & & & & & & & & & & & & \\
            
            \multirow{3}[0]{*}{200-299}
            & A &  & 0.967 & 0.381 & 248.55 &  & 0.423 & 0.174 & 248.55 &  & \textbf{0.025} & 0.005 & 248.55 \\ 
            & B &  & 0.904 & 0.386 & 248.55 &  & 0.421 & 0.165 & 248.55 &  & \textbf{0.034} & 0.009 & 248.55 \\ 
            & C &  & 1.181 & 0.429 & 248.55 &  & 0.438 & 0.166 & 248.55 &  & \textbf{0.067} & 0.008 & 248.55 \\ 
            & & & & & & & & & & & & \\
         
            \multirow{3}[0]{*}{300-399}
            & A &  & 1.245 & 0.594 & 341.67 &  & 0.513 & 0.166 & 341.67 &  & \textbf{0.045} & 0.010 & 341.67 \\ 
            & B &  & 1.531 & 0.671 & 341.67 &  & 0.617 & 0.321 & 341.67 &  & \textbf{0.171} & 0.057 & 341.67 \\ 
            & C &  & 1.954 & 1.023 & 341.67 &  & 0.989 & 0.348 & 341.67 &  & \textbf{0.370} & 0.151 & 341.67 \\ 
            & & & & & & & & & & \\
            
            \multirow{3}[0]{*}{400-499}
            & A &  & 2.126 & 1.292 & 444.60 &  & 1.172 & 0.495 & 444.60 &  & \textbf{0.140} & 0.021 & 444.60 \\ 
            & B &  & 2.035 & 0.852 & 444.60 &  & 0.997 & 0.479 & 444.60 &  & \textbf{0.229} & 0.058 & 444.60 \\ 
            & C &  & 2.521 & 1.267 & 444.60 &  & 1.626 & 0.855 & 444.60 &  & \textbf{0.619} & 0.180 & 444.60 \\ 
            & & & & & & & & & & \\
            
            \multirow{3}[0]{*}{500-599}
            & A &  & 2.585 & 1.764 & 548.56 &  & 1.447 & 0.931 & 548.56 &  & \textbf{0.517} & 0.183 & 548.56 \\ 
            & B &  & 2.790 & 1.752 & 548.56 &  & 1.323 & 0.519 & 548.57 &  & \textbf{0.605} & 0.241 & 548.56 \\ 
            & C &  & 2.966 & 1.930 & 548.56 &  & 1.603 & 0.934 & 548.56 &  & \textbf{0.880} & 0.200 & 548.56 \\ 
            & & & & & & & & & & \\
            
            \multirow{3}[0]{*}{600-699}
            & A &  & 2.667 & 1.657 & 648.34 &  & 1.351 & 0.671 & 648.33 &  & \textbf{0.443} & 0.150 & 648.33 \\ 
            & B &  & 2.922 & 1.929 & 648.34 &  & 1.496 & 0.781 & 648.33 &  & \textbf{0.657} & 0.224 & 648.33 \\ 
            & C &  & 3.277 & 2.358 & 648.34 &  & 1.859 & 0.896 & 648.33 &  & \textbf{1.294} & 0.663 & 648.33 \\ 
            & & & & & & & & & & \\
            
            \multirow{3}[0]{*}{700-799} 
            & A &  & 2.765 & 1.720 & 741.01 &  & 1.404 & 0.805 & 741.00 &  & \textbf{0.874} & 0.311 & 741.00 \\ 
            & B &  & 3.259 & 2.317 & 741.01 &  & 1.800 & 0.904 & 741.02 &  & \textbf{1.005} & 0.534 & 741.00 \\ 
            & C &  & 3.261 & 1.847 & 741.01 &  & 1.667 & 0.827 & 741.00 &  & \textbf{1.342} & 0.643 & 741.00 \\ 
            & & & & & & & & & & \\
            
            \multirow{3}[0]{*}{800-899}
            & A &  & 3.045 & 2.158 & 847.01 &  & 1.530 & 0.978 & 847.00 &  & \textbf{0.893} & 0.253 & 847.00 \\ 
            & B &  & 3.365 & 2.263 & 847.01 &  & 1.755 & 0.896 & 847.00 &  & \textbf{1.173} & 0.483 & 847.00 \\ 
            & C &  & 3.454 & 2.325 & 847.01 &  & \textbf{1.680} & 0.472 & 847.00 &  & 1.899 & 1.100 & 847.00 \\ 
            & & & & & & & & & & \\
            
            \multirow{3}[0]{*}{900-1001} 
            & A &  & 3.268 & 2.258 & 957.41 &  & 1.359 & 0.541 & 957.40 &  & \textbf{1.012} & 0.349 & 957.40 \\ 
            & B &  & 3.284 & 2.282 & 957.42 &  & 1.363 & 0.563 & 957.40 &  & \textbf{1.246} & 0.623 & 957.40 \\ 
            & C &  & 2.616 & 1.161 & 957.42 &  & \textbf{0.991} & 0.065 & 957.40 &  & 2.016 & 1.022 & 957.40 \\ 
            & & & & & & & & & & \\
            
            ALL 
            & A &  & 1.569 & 0.934 & 412.80 &  & 0.777 & 0.386 & 412.80 &  & \textbf{0.256} & 0.081 & 412.80 \\ 
            & B &  & 1.655 & 0.942 & 412.80 &  & 0.788 & 0.364 & 412.80 &  & \textbf{0.343} & 0.144 & 412.80 \\ 
            & C &  & 1.860 & 1.031 & 412.80 &  & 0.944 & 0.413 & 412.80 &  & \textbf{0.587} & 0.257 & 412.80 \\ 
            \midrule
			\multicolumn{2}{c}{Average}  &  & 1.695 & 0.969 & 412.80 &  & 0.837 & 0.388 & 412.80 &  & \textbf{0.395} & 0.161 & 412.80 \\ 
			\bottomrule
		\end{tabular}
	}
\end{table}

From these results, we first observe that our re-implementation of the RR algorithm has found solutions of a similar or better quality than the original RR-VRVC. As seen in Table~\ref{tab:summaryres1}, for all instances with 200 or more services, RR finds solutions of a better quality than RR-VRVC in a shorter amount of time. Generally, the use of the termination criterion and cooling schedule proportional to the size of the instances is beneficial for the method.

The RR solution approach is further improved by relying on the fast evaluation strategy for the pickup-and-delivery insertions (RR-fast). The increased speed of these evaluations permits to do significantly more iterations, leading to an average gap of $-0.085\%$ for RR-fast compared to $0.164\%$ for RR on the classic instances, and $0.837\%$ compared to $1.695\%$ on the X instances. The significance of these solution-quality improvement is confirmed by pairwise Wilcoxon tests, with p-values of $4.90 \times 10^{-60}$ and $2.06 \times 10^{-286}$ on the classic and X instance sets, respectively. Notably, this means that the proposed mechanisms for evaluating insertions are not only useful from a complexity viewpoint, but also from a practical perspective.

Overall, considering all methods, the HGS with the proposed local search achieves the best performance with an average gap of $-0.290\%$ on the classical instances, and $0.395\%$ on the X instances. The significance of these improvements over RR-fast is also confirmed by pairwise Wilcoxon tests, with respective p-values of $4.24 \times 10^{-67}$ and $2.67 \times 10^{-289}$. Out of the 128 instances from \cite{Renaud2000a}, it found or improved the best-known solution in 126 cases (i.e., all BKS except those of KROD99C and PR107B, reported only in RBO). Moreover, among these 126 cases, the BKS has been matched on all ten runs in 82 cases, whereas the average solution quality of HGS was better than the BKS for the remaining 44 cases. Finally, the 28 known optimal solutions from \cite{Dumitrescu2010} have been found on all runs.

Solution-quality differences are more marked for larger instances with many services. For example, an improvement of $3.768\%$ over previous known results have been achieved for instance PR439C. HGS produces superior results for the wide majority of the instances, though on a few of the largest instances of set X (class C with more than 800 visits), the simpler RR-fast strategy seems to perform better. This is likely due to the fact that, on such instances, the crossover mechanism of HGS creates many disruptions in the new solutions, which are not fully re-optimized in a single local-search descent. In contrast, the ruin-and-recreate mechanism is more localized. This effect is also exacerbated by the precedence constraints implied by the pickups and deliveries, which limit the amount of feasible moves.

\subsection{Impact of the Neighborhoods}
\label{sec:exp:sensitivity}

To understand the contribution of each of the six neighborhoods used in the local search of HGS, we conduct an ablation study which consists of removing one neighborhood at a time and comparing the performance of the resulting approach to our baseline results, using the same benchmark instances and stopping criterion as in Section~\ref{tab:lscalibration}. Table~\ref{tab:ablation-study} displays summarized results for this experiment over the instances of set X. It indicates, for each algorithm variant, the average percentage gap over ten runs relative to the BKS, as well as the p-value of a pairwise Wilcoxon test between the considered configuration and the baseline without any ablation.

\begin{table}[htbp]
	\centering
	\caption{Ablation study on the impact of the different neighborhoods}
	\label{tab:ablation-study}
    \setlength\tabcolsep{7pt}
    \renewcommand{\arraystretch}{1.2}
	\scalebox{0.85}
	{
	\begin{tabular}{|l|c|c|}
		\hline
	    Configuration & Average Gap(\%) & p-value \\
		\hline
        Baseline & \textbf{0.395} & --- \\ 
        No \textsc{Relocate Pair}  & 0.463 & \num{5.63703274485541E-21} \\ 
        No \textsc{2-Opt}  & 0.399 &  \num{0.891003707448502} \\ 
        No \textsc{2k-Opt}  & 0.405 & \num{0.600254641297334} \\ 
        No \textsc{Or-Opt} & 0.793 & \num{2.83671839473107E-216} \\ 
        No \textsc{4-Opt} & 0.409 & \num{0.072601111024658}  \\ 
        No \textsc{Balas-Simonetti} & 0.416 & \num{0.051551455963718} \\ 
        No Large Neighborhoods & 0.423 & \num{0.003229077825285} \\ 
		\hline
	\end{tabular}%
	}
\end{table}%

As seen in this experiment, the baseline method with all the neighborhoods achieved better average gaps than the other configurations. In terms of their individual impact, canonical neighborhoods such as \textsc{Or-Opt} and \textsc{Relocate Pair} are indispensable to achieve good solutions. This is confirmed by pairwise Wilcoxon tests comparing the results of the configurations obtained by deactivating these neighborhoods with those of the baseline.

The \textsc{Relocate Pair}, \textsc{2-Opt}, and  \textsc{Or-Opt} neighborhoods are explored in the first stage in our local search (Lines 3 to 6 of Algorithm~\ref{Algo:LS}) as it is possible to decompose their exploration and rapidly apply improving moves. In contrast, the large neighborhoods are explored less frequently since they require a complete quadratic-time exploration. Their main goal is to help escaping local minima and sub-optimal structures that may otherwise persist in the solutions. These distinct roles are clearly visible as we monitor the local search. As seen in Table~EC.5 (in the Electronic Companion),
a local search on an instance of Class X applies on average $21.48$ \textsc{Relocate Pair}, $2.95$ \textsc{2-Opt}, and $16.43$ \textsc{Or-Opt} moves, compared to only $1.71$ \textsc{2k-Opt}, $1.86$ \textsc{4-Opt} and $1.45$ \textsc{Balas-Simonetti} moves. Still, a Wilcoxon test demonstrates that jointly deactivating all the large neighborhoods leads to results that are significantly worse. These observations are also in line with other studies on variable neighborhood search and large neighborhood search, which have shown that a diversity neighborhood structures permits to effectively escape low-quality local minima, since a local minimum of a specific neighborhood may not be a minimum of other ones \citep{Hansen2019,Pisinger2019}.

\subsection{Speed evaluation on Grubhub data}
\label{sec:exp:mealdelivery}

Finally, we analyze the performance of the complete HGS method on the Grubhub data set, which contains instances with $n \in \{2,\dots,15\}$ pickup-and-delivery pairs and therefore $N \in \{5,\dots,31\}$ visits. These instances involve open tours (i.e., the return to the origin point is not counted). There are ten instances for each size, totaling 140 instances. Solution time is critical since these instances relate to a real-time meal delivery application. Ideally, it should be no greater than a fraction of a second. For this reason, we use a short termination criterion of 100 successive iterations without improvement.

\cite{ONeil2018a} and \cite{ONeil2019} have proposed different solution techniques to solve these instances to near optimality. The authors suggested using constraint programming, decision diagrams, and integer programming to achieve optimality and terminate the search quickly. Based on their experiments, they proposed a configuration called ``Ruland+CP+AP'', which consists of using constraint programming as a hot start for a MIP solver based on the formulation from \cite{Ruland1997}. They conducted experiments on an Intel Core i7-7300U @2.6GHz CPU and measured the time taken to 
(i) attain a solution within 10\% of the optimum, 
(ii) attain a solution within 5\% of the optimum,
(iii) attain the optimal solution,
(iv) prove optimality.

\begin{figure}[htbp]
\centering
\includegraphics[width=\linewidth]{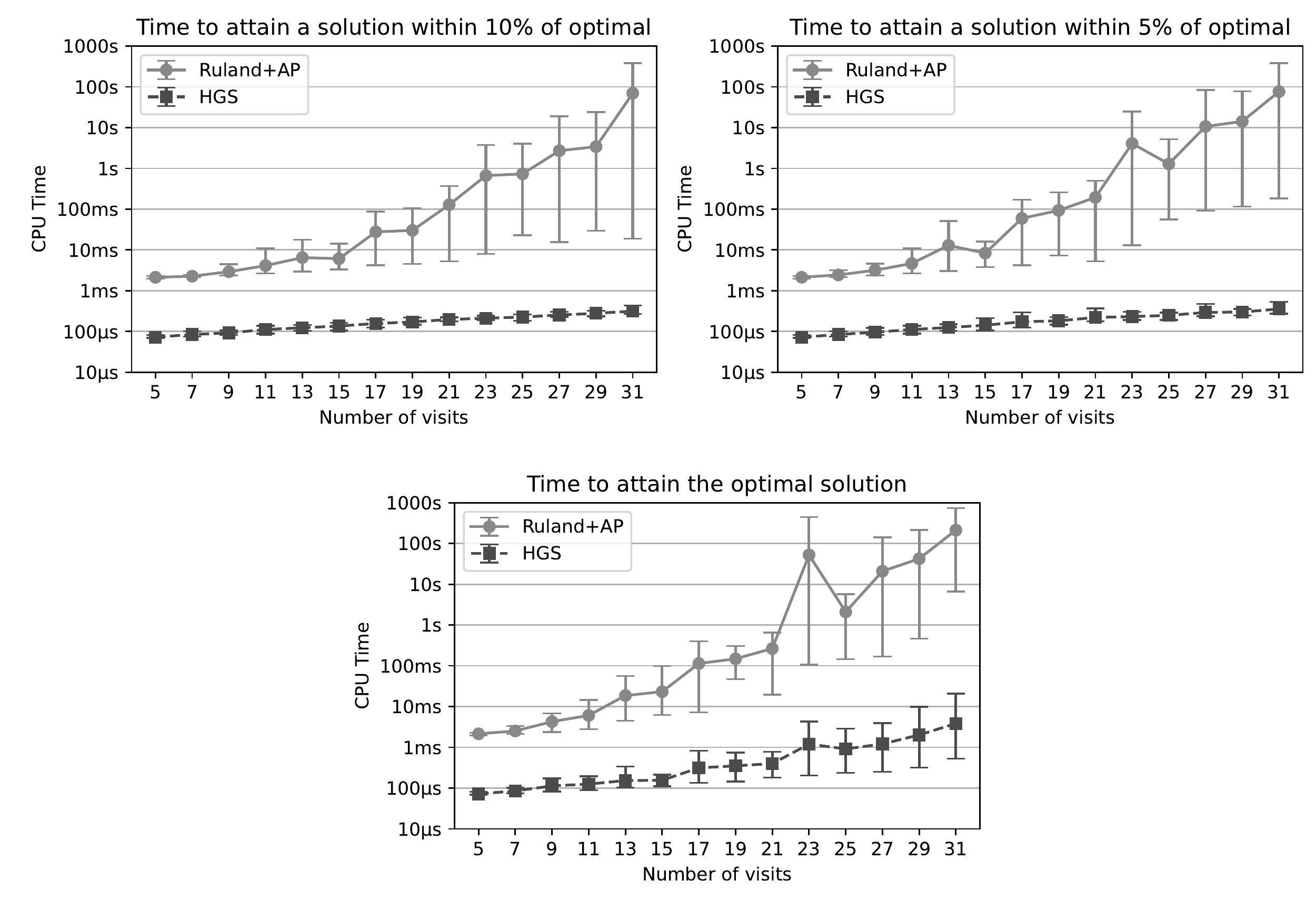}
\caption{Minimum, average and maximum computational effort to attain solutions within 10\%, 5\% and 0\% of the optimum}
\label{fig:grubseries}
\end{figure}

\begin{figure}[htbp]
\centering
\includegraphics[width=0.8\linewidth]{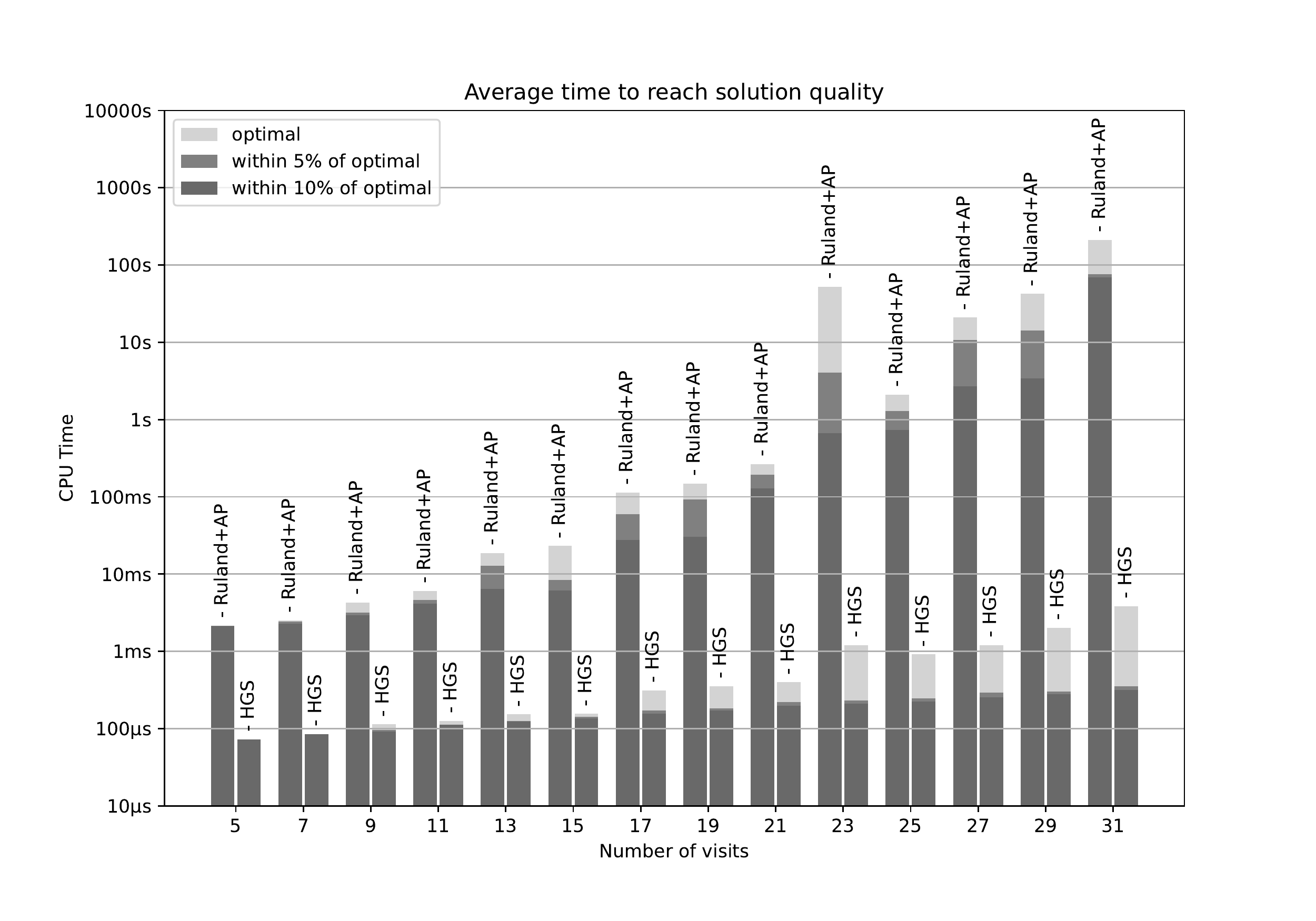}
\caption{Distribution of the CPU time needed to locate solutions within 10\%, 5\% and 0\% of the optimum}
\label{fig:grubstack}
\end{figure}

As we rely on heuristics, we do not prove optimality, but we can nonetheless measure in a similar fashion the time to achieve (i), (ii), and (iii). Figures~\ref{fig:grubseries}~and~\ref{fig:grubstack} compare the CPU time of HGS and Ruland+CP+AP to achieve those targets. For consistency, since Ruland+CP+AP performed a single run per instance, we used the first run of HGS to create the figures. The CPU time of Ruland+CP+AP was also scaled by a factor $1.086$ to account for CPU speed differences. For each group of ten instances of each given size, the first figure reports the minimum, average, and maximum time needed to attain each target. The second figure depicts the average solution time needed to attain each target as a stacked bar plot. In both cases, CPU times are represented on a logarithmic scale. Finally, Tables~EC.6~and~EC.7
(in the Electronic Companion) provide the detailed time and solution values for each instance, considering the ten runs of HGS.

As observed in these experiments, the proposed HGS with local search attained the known optimal solution for all instances and runs (i.e., on 1400 executions overall). This confirms the performance and robustness of the approach. Most notably, solutions within 5\% of the optimum are typically found within a fraction of a millisecond, and optimal solutions are found in milliseconds, for all instances, as visible in Figure~\ref{fig:grubseries}. This is a remarkable capability, given that instances with $15$ pickup-and-delivery services (i.e., $31$ visits) are far from trivial, with already $30!/2^{15} \approx 8.09 \times 10^{27}$ possible solutions. In comparison, the Ruland+CP+AP approach based on constraint programming and integer programming requires dozens of seconds on average to achieve the same targets.

We also measured the number of local search runs needed to attain the same targets. Interestingly, we observed that a single local search run was sufficient to attain the optimal solution in 61.79\% of the cases. Similarly, 99.93\% of the optimal solutions were found when generating the initial population, whereas the remaining 0.07\% solutions were found during the evolution of the genetic algorithm on these instances. This indicates that the proposed local search effectively locates high-quality PDTSP solutions and that it represents a meaningful alternative for applications requiring very short response time.

\section{Conclusions} 
\label{chap:6}

In this paper, we have introduced efficient search strategies and new large neighborhoods for the pickup-and-delivery traveling salesman problem. Our study permitted significant complexity reductions for classical neighborhoods such as \textsc{Relocate Pair} as well as the efficient exploration of large neighborhoods such as \textsc{2k-Opt}, \textsc{4-Opt} and \textsc{Balas-Simonetti}. We conducted extensive computational experiments, demonstrating that the use of these neighborhoods in the classical HGS metaheuristic leads to a robust search algorithm that produces almost systematically the best-known and optimal solutions for all the considered instances. Our methodology also permitted significant speed-ups and improvements over a ruin-and-recreate algorithm previously proposed by \cite{Veenstra2017}. Our search approach is especially useful in applications where time is critical, e.g., meal-delivery and mobility-on-demand, as it permits us to repeatedly achieve optimal or near-optimal solutions for small problems within milliseconds.

The research perspectives are numerous. First of all, research can be pursued on some of the neighborhood structures discussed in this paper, such as \textsc{2k-Opt} and \textsc{4-Opt}. Since the evaluation of these neighborhoods is grounded on dynamic programming, identifying the families of problem attributes \citep{Vidal2012a} for which an efficient (sub-quadratic) search is achievable is an important and challenging research question. Studies on incremental solution strategies and decompositions for the dynamic programs could also lead to substantial improvements and give more flexibility on move acceptance policies.

More generally, research should be pursued on fundamental neighborhood structures for canonical problems in the vehicle routing family. Many studies have been focused on high-level search strategies (often new metaheuristic concepts) over the past two decades, but improvements along this line see diminishing returns nowadays. The core component of state-of-the-art metaheuristics is often a low-complexity neighborhood search, used for systematic solution improvement. However, the best possible complexity for many useful neighborhoods still remains an open research question, and regular breakthroughs are still being made \citep{DeBerg2020}. Generally, we believe that a  disciplined study of neighborhood searches, grounded on computational complexity theory, combinatorial optimization, and possibly machine learning techniques \citep{Arnold2021,Karimi-Mamaghan2021,Lancia2019}, still has potential for significant breakthrough and paradigm shifts.

\ACKNOWLEDGMENT{The authors would like to thank Marjoleen Aerts-Veenstra, Jacques Renaud, and Ryan J. O'Neil for their prompt answers and detailed results, which were instrumental to complete the numerical comparisons of this paper. This research has been partially funded by CAPES, CNPq [grant number 308528/2018-2], and FAPERJ [grant number E-26/202.790/2019] in Brazil. This research was also enabled in part by computational resources provided by Calcul Qu\'ebec and Compute Canada. This support is gratefully acknowledged.
}

\ECSwitch

\ECHead{Electronic Companion: Exponential-Size Neighborhoods for the Pickup-and-Delivery Traveling Salesman Problem}

This Electronic Companion provides additional detailed results. The instances of Class~1~and~2 have been introduced in \cite{Renaud2000a}, whereas the instances of Class 3 come from \cite{Dumitrescu2010}.
The Grubhub data sets have been introduced in \cite{ONeil2018a} and \cite{ONeil2019} and are accessible at \url{https://github.com/grubhub/tsppdlib}. Finally, the instances of Class X have been introduced in the present study.

Tables \ref{tab:compresX}--\ref{tab:compres3} and \ref{tab:compresgrubhub} provide detailed results over ten runs for all classes of instances. The format of the tables is identical to that of Table~5
in the main paper, except that Best and Avg solution values as well as CPU times (T(s)) are provided for each separate instance. Known optimal solutions from \cite{Dumitrescu2010} are also marked with an~$^*$.
As discussed in the paper, CPU times have been scaled to obtain an equivalent time on an AMD Rome 7532 2.40GHz.
Table~\ref{tab:summaryresNBMoves} reports the average number of moves of each type applied during a local search on the different instances subgroups.
Finally, Table~\ref{tab:compresgrubhubtimes} gives the detailed data used to generate Figures~7~and~8,
based on the the first run of HGS. It compares the time taken by Ruland+CP+AP and HGS to attain a solution-quality target of 10\%, 5\%, and 0\% from the known optimal solution.

\begin{landscape}

{
	\tiny
	\renewcommand{\arraystretch}{0.99}
\begin{longtable}{@{}lcccccccccccccc@{}}
		\caption{Detailed results on the instances of Class~X} \label{tab:compresX}\\
		\toprule
		& & & & \multicolumn{3}{c}{RR} & & \multicolumn{3}{c}{RR-FAST} & & \multicolumn{3}{c}{HGS}\\
		\cmidrule{5-7}\cmidrule{9-11}\cmidrule{13-15}
		Instance & N & BKS &
		& Avg & Best & T(s) & 
		& Avg & Best & T(s) & 
		& Avg & Best & T(s) \\  
		\midrule
		\endhead
		\hline \multicolumn{15}{r}{\textit{(continued next page)}} \\
		\endfoot
		\bottomrule
		\endlastfoot
		X-n101-k25A & 101 & 8333 &  & \textbf{8333.00} & 8333.00 & 101.00 &  & \textbf{8333.00} & 8333.00 & 101.00 &  & \textbf{8333.00} & 8333.00 & 101.00 \\ 
        X-n101-k25B & 101 & 8341 &  & \textbf{8341.00} & 8341.00 & 101.00 &  & \textbf{8341.00} & 8341.00 & 101.00 &  & \textbf{8341.00} & 8341.00 & 101.00 \\ 
        X-n101-k25C & 101 & 8549 &  & \textbf{8549.00} & 8549.00 & 101.00 &  & \textbf{8549.00} & 8549.00 & 101.00 &  & \textbf{8549.00} & 8549.00 & 101.00 \\ 
        X-n106-k14A & 106 & 6019 &  & \textbf{6019.00} & 6019.00 & 105.00 &  & \textbf{6019.00} & 6019.00 & 105.00 &  & \textbf{6019.00} & 6019.00 & 105.00 \\ 
        X-n106-k14B & 106 & 6288 &  & \textbf{6288.00} & 6288.00 & 105.00 &  & \textbf{6288.00} & 6288.00 & 105.00 &  & \textbf{6288.00} & 6288.00 & 105.00 \\ 
        X-n106-k14C & 106 & 6399 &  & \textbf{6399.00} & 6399.00 & 105.00 &  & \textbf{6399.00} & 6399.00 & 105.00 &  & \textbf{6399.00} & 6399.00 & 105.00 \\ 
        X-n110-k13A & 110 & 9421 &  & \textbf{9421.00} & 9421.00 & 109.00 &  & \textbf{9421.00} & 9421.00 & 109.00 &  & \textbf{9421.00} & 9421.00 & 109.00 \\ 
        X-n110-k13B & 110 & 9429 &  & \textbf{9429.00} & 9429.00 & 109.00 &  & \textbf{9429.00} & 9429.00 & 109.00 &  & \textbf{9429.00} & 9429.00 & 109.00 \\
\end{longtable}
}%

\pagebreak
{
	\tiny
	\begin{longtable}{@{}lcccccccccccccccccccc@{}}
		\caption{Detailed Results on the instances of Class~1} \label{tab:compres1}\\
		\toprule
		& & & & \multicolumn{1}{c}{RBO/BC10} & & \multicolumn{3}{c}{RR-VRVC} & & \multicolumn{3}{c}{RR} & & \multicolumn{3}{c}{RR-FAST} & & \multicolumn{3}{c}{HGS}\\
		\cmidrule{5-5}\cmidrule{7-9}\cmidrule{11-13}\cmidrule{15-17}\cmidrule{19-21}
		Instance & N & BKS & & Best  &
		& Avg & Best & T(s) & 
		& Avg & Best & T(s) & 
		& Avg & Best & T(s) & 
		& Avg & Best & T(s) \\  
		\midrule
		\endhead
		\hline \multicolumn{21}{r}{\textit{(continued next page)}} \\
		\endfoot
		\bottomrule
		\endlastfoot
        EIL51A & 51 & 464* &  & 464 &  & \textbf{464.00} & 464.00 & 3.62 &  & \textbf{464.00} & 464.00 & 51.00 &  & \textbf{464.00} & 464.00 & 51.00 &  & \textbf{464.00} & 464.00 & 51.00 \\ 
        EIL51B & 51 & 469* &  & 469 &  & \textbf{469.00} & 469.00 & 3.83 &  & \textbf{469.00} & 469.00 & 51.00 &  & \textbf{469.00} & 469.00 & 51.00 &  & \textbf{469.00} & 469.00 & 51.00 \\ 
        EIL51C & 51 & 488 &  & 488 &  & \textbf{488.00} & 488.00 & 3.90 &  & \textbf{488.00} & 488.00 & 51.00 &  & \textbf{488.00} & 488.00 & 51.00 &  & \textbf{488.00} & 488.00 & 51.00 \\ 
        ST69A & 69 & 764* &  & 764 &  & \textbf{764.00} & 764.00 & 8.16 &  & \textbf{764.00} & 764.00 & 69.00 &  & \textbf{764.00} & 764.00 & 69.00 &  & \textbf{764.00} & 764.00 & 69.00 \\ 
        ST69B & 69 & 771 &  & 771 &  & \textbf{771.00} & 771.00 & 8.09 &  & \textbf{771.00} & 771.00 & 69.00 &  & \textbf{771.00} & 771.00 & 69.00 &  & \textbf{771.00} & 771.00 & 69.00 \\ 
        ST69C & 69 & 793 &  & 793 &  & \textbf{793.00} & 793.00 & 9.08 &  & \textbf{793.00} & 793.00 & 69.00 &  & \textbf{793.00} & 793.00 & 69.00 &  & \textbf{793.00} & 793.00 & 69.00 \\ 
        EIL75A & 75 & 583 &  & 583 &  & 589.30 & 583.00 & 10.50 &  & \textbf{583.00} & 583.00 & 75.00 &  & \textbf{583.00} & 583.00 & 75.00 &  & \textbf{583.00} & 583.00 & 75.00 \\ 
        EIL75B & 75 & 601 &  & 601 &  & 601.30 & 601.00 & 11.28 &  & \textbf{601.00} & 601.00 & 75.00 &  & \textbf{601.00} & 601.00 & 75.00 &  & \textbf{601.00} & 601.00 & 75.00 \\
\end{longtable}
}%

\pagebreak
\begin{table}[htbp]
	\scriptsize
	\centering
	\caption{Detailed results on the instances of Class~2}
	\label{tab:compres2}
	~\\[2mm]
	\begin{tabular}{@{}lcccccccccccccccccccc@{}}
		\toprule
		& & & & \multicolumn{1}{c}{RBO/BC10} & & \multicolumn{3}{c}{RR-VRVC} & & \multicolumn{3}{c}{RR} & & \multicolumn{3}{c}{RR-FAST} & & \multicolumn{3}{c}{HGS}\\
		\cmidrule{5-5}\cmidrule{7-9}\cmidrule{11-13}\cmidrule{15-17}\cmidrule{19-21}
		Instance & N & BKS & & Best &
		& Avg & Best & T(s) & 
		& Avg & Best & T(s) & 
		& Avg & Best & T(s) & 
		& Avg & Best & T(s) \\  
		\midrule
        N101P1 & 101 & 799* &  & \textbf{799} &  & \textbf{799.00} & 799.00 & 26.38 &  & \textbf{799.00} & 799.00 & 101.00 &  & \textbf{799.00} & 799.00 & 101.00 &  & \textbf{799.00} & 799.00 & 101.00 \\ 
        N101P2 & 101 & 729* &  & \textbf{729} &  & \textbf{729.00} & 729.00 & 26.95 &  & \textbf{729.00} & 729.00 & 101.00 &  & \textbf{729.00} & 729.00 & 101.00 &  & \textbf{729.00} & 729.00 & 101.00 \\ 
        N101P3 & 101 & 748* &  & \textbf{748} &  & \textbf{748.00} & 748.00 & 26.88 &  & \textbf{748.00} & 748.00 & 101.00 &  & \textbf{748.00} & 748.00 & 101.00 &  & \textbf{748.00} & 748.00 & 101.00 \\ 
        N101P4 & 101 & 807* &  & \textbf{807} &  & \textbf{807.00} & 807.00 & 26.88 &  & \textbf{807.00} & 807.00 & 101.00 &  & \textbf{807.00} & 807.00 & 101.00 &  & \textbf{807.00} & 807.00 & 101.00 \\ 
        N101P5 & 101 & 783* &  & \textbf{783} &  & \textbf{783.00} & 783.00 & 26.88 &  & \textbf{783.00} & 783.00 & 101.00 &  & \textbf{783.00} & 783.00 & 101.00 &  & \textbf{783.00} & 783.00 & 101.00 \\ 
        N101P6 & 101 & 755* &  & \textbf{755} &  & \textbf{755.00} & 755.00 & 26.95 &  & \textbf{755.00} & 755.00 & 101.00 &  & \textbf{755.00} & 755.00 & 101.00 &  & \textbf{755.00} & 755.00 & 101.00 \\
\bottomrule
	    \end{tabular}%
\end{table}%

\pagebreak
\begin{table}[htbp]
	\scriptsize
	\centering
	\caption{Detailed results on the instances of Class 3}
	\label{tab:compres3}
	~\\[2mm]
	\begin{tabular}{@{}lclcccccccccccccccccc@{}}
		\toprule
		& & & & \multicolumn{1}{c}{BC10} & & \multicolumn{3}{c}{RR-VRVC} & & \multicolumn{3}{c}{RR} & & \multicolumn{3}{c}{RR-FAST} & & \multicolumn{3}{c}{HGS}\\
		\cmidrule{5-5}\cmidrule{7-9}\cmidrule{11-13}\cmidrule{15-17}\cmidrule{19-21}
		Instance & N & BKS & & Best &
		& Avg & Best & T(s) & 
		& Avg & Best & T(s) & 
		& Avg & Best & T(s) & 
		& Avg & Best & T(s) \\  
		\midrule
        prob5a & 11 & 3585* &  & 3585 &  & \textbf{3585.00} & 3585.00 & 0.21 &  & \textbf{3585.00} & 3585.00 & 11.00 &  & \textbf{3585.00} & 3585.00 & 11.00 &  & \textbf{3585.00} & 3585.00 & 11.00 \\ 
        prob5b & 11 & 2565* &  & 2565 &  & \textbf{2565.00} & 2565.00 & 0.07 &  & 2807.00 & 2807.00 & 11.00 &  & 2807.00 & 2807.00 & 11.00 &  & \textbf{2565.00} & 2565.00 & 11.00 \\ 
        prob5c & 11 & 3787* &  & 3787 &  & \textbf{3787.00} & 3787.00 & 0.07 &  & \textbf{3787.00} & 3787.00 & 11.00 &  & \textbf{3787.00} & 3787.00 & 11.00 &  & \textbf{3787.00} & 3787.00 & 11.00 \\ 
        prob5d & 11 & 3128* &  & 3128 &  & \textbf{3128.00} & 3128.00 & 0.07 &  & 3230.00 & 3230.00 & 11.00 &  & 3230.00 & 3230.00 & 11.00 &  & \textbf{3128.00} & 3128.00 & 11.00 \\ 
        prob5e & 11 & 3123* &  & 3123 &  & \textbf{3123.00} & 3123.00 & 0.07 &  & \textbf{3123.00} & 3123.00 & 11.00 &  & \textbf{3123.00} & 3123.00 & 11.00 &  & \textbf{3123.00} & 3123.00 & 11.00 \\ 
        prob10a & 21 & 4896* &  & 4896 &  & \textbf{4896.00} & 4896.00 & 0.35 &  & \textbf{4896.00} & 4896.00 & 21.00 &  & \textbf{4896.00} & 4896.00 & 21.00 &  & \textbf{4896.00} & 4896.00 & 21.00 \\
\bottomrule
	\end{tabular}%
\end{table}%

\end{landscape}

\vspace*{0cm}
\begin{table}[htbp]
	\centering
	\caption{Average number of moves per local search on the instances of Class X} \label{tab:summaryresNBMoves}
	\setlength\tabcolsep{6.5pt}
	\renewcommand{\arraystretch}{1.05}
	\scalebox{0.8}
	{
		\begin{tabular}{ccC{1.6cm}C{1.6cm}C{1.6cm}C{1.6cm}C{1.6cm}C{1.6cm}}
			\toprule
			$N$ & Type & \textsc{Reloc} & \textsc{2-Opt} & \textsc{Or-Opt} & \textsc{2k-Opt} & \textsc{4-Opt} & \textsc{B\&S} \\
			\midrule
			\multirow{3}[0]{*}{100-199}
             & A &  9.02 & 1.34 & 8.23 & 1.42 & 1.64 & 1.30 \\
             & B &  11.51 & 1.77 & 9.90 & 1.52 & 1.66 & 1.34 \\
             & C &  15.66 & 2.75 & 12.68 & 1.64 & 1.70 & 1.40 \\
             & & & & & & & \\
             \multirow{3}[0]{*}{200-299}
             & A & 14.54 & 1.88 & 12.16 & 1.56 & 1.83 & 1.39 \\
             & B & 18.36 & 2.45 & 14.37 & 1.67 & 1.82 & 1.44 \\
             & C &  27.33 & 4.23 & 19.73 & 1.88 & 1.91 & 1.54 \\
\bottomrule
		\end{tabular}
 }
\end{table}

\pagebreak

{
	\scriptsize
	\begin{longtable}{@{}lcccccc@{}}
		\caption{Detailed results on the Grubhub instances} \label{tab:compresgrubhub}\\
		\toprule
		& & & &\multicolumn{3}{c}{HGS} \\
		\cmidrule{5-7}
		Instance & N & BKS & & Avg & Best & T(s) \\
		\midrule
		\endhead
		\hline 
		\multicolumn{7}{r}{\textit{(continued on next page)}} \\
		\endfoot
		\bottomrule
		\endlastfoot
        grubhub-02-0 & 5 & 3214* & & \textbf{3214} & 3214 & \num{0.0044853} \\
        grubhub-02-1 & 5 & 2736* & & \textbf{2736} & 2736 & \num{0.0045313} \\
        grubhub-02-2 & 5 & 2558* & & \textbf{2558} & 2558 & \num{0.0044867} \\
        grubhub-02-3 & 5 & 1382* & & \textbf{1382} & 1382 & \num{0.0044228} \\
        grubhub-02-4 & 5 & 2467* & & \textbf{2467} & 2467 & \num{0.0043727} \\
        grubhub-02-5 & 5 & 3099* & & \textbf{3099} & 3099 & \num{0.0044626} \\
        grubhub-02-6 & 5 & 1778* & & \textbf{1778} & 1778 & \num{0.0043754} \\
        grubhub-02-7 & 5 & 2206* & & \textbf{2206} & 2206 & \num{0.0044422} \\
\end{longtable}
}%

{
	\scriptsize
	\begin{longtable}{cccccccccc}
		\caption{Time-to-target results on Grubhub instances} \label{tab:compresgrubhubtimes}\\
		\toprule
		& & &\multicolumn{3}{c}{T(s) of Ruland+AP} & &\multicolumn{3}{c}{T(s) of HGS} \\
		\cmidrule{4-6} \cmidrule{8-10}
		Instance & N & & 10\% & 5\% & 0\% & & 10\% & 5\% & 0\% \\
		\midrule
		\endhead
		\hline 
		\multicolumn{10}{r}{\textit{(continued on next page)}} \\
		\endfoot
		\bottomrule
		\endlastfoot
        grubhub-02-0 & 5 &  & \num{0.0022771914} & \num{0.0022771914} & \num{0.0022771914} &  & \num{0.0000728000} & \num{0.0000728000} & \num{0.0000728000} \\ 
        grubhub-02-1 & 5 &  & \num{0.0020829766} & \num{0.0020829766} & \num{0.0022063168} &  & \num{0.0000695000} & \num{0.0000695000} & \num{0.0000695000} \\ 
        grubhub-02-2 & 5 &  & \num{0.0022127599} & \num{0.0022127599} & \num{0.0022127599} &  & \num{0.0000723000} & \num{0.0000723000} & \num{0.0000723000} \\ 
        grubhub-02-3 & 5 &  & \num{0.0020774539} & \num{0.0021925101} & \num{0.0021925101} &  & \num{0.0000716000} & \num{0.0000716000} & \num{0.0000716000} \\ 
        grubhub-02-4 & 5 &  & \num{0.0022219644} & \num{0.0022219644} & \num{0.0022219644} &  & \num{0.0000695000} & \num{0.0000695000} & \num{0.0000695000} \\ 
        grubhub-02-5 & 5 &  & \num{0.0021041469} & \num{0.0021041469} & \num{0.0021041469} &  & \num{0.0000694000} & \num{0.0000694000} & \num{0.0000694000} \\ 
        grubhub-02-6 & 5 &  & \num{0.0021455672} & \num{0.0021455672} & \num{0.0021455672} &  & \num{0.0000724000} & \num{0.0000724000} & \num{0.0000724000} \\ 
        grubhub-02-7 & 5 &  & \num{0.0019688409} & \num{0.0019688409} & \num{0.0019688409} &  & \num{0.0000814000} & \num{0.0000814000} & \num{0.0000814000} \\
\end{longtable}
}%

\end{document}